\documentclass[preprint, 11pt,authoryear]{LaTeXConfig/elsarticle} %

\usepackage[utf8]{inputenc}
\usepackage[british]{babel}

\usepackage{tikz}
\usepackage{pgfplots}
\pgfplotsset{compat=1.18}
\usepackage{relsize}
\usepackage{graphicx}
\usepackage{caption}
\usepackage{subcaption}
\usepackage{enumitem}
\usepackage{xcolor}

\usepackage{amsmath,amsfonts,amssymb, mathtools}
\mathtoolsset{showonlyrefs=true,showmanualtags=true} %
\usepackage{amsthm}

\newtheorem{definition}{Definition}[section]
\newtheorem{corollary}{Corollary}[section]

\newtheorem{prop}{Proposition}[section]

\usepackage[linesnumbered,vlined,ruled,commentsnumbered]{algorithm2e}

\usepackage{booktabs}
\usepackage{tabularx}

\usepackage{numprint}

\usepackage[short]{optidef}

\usepackage{listings}

\usepackage[hidelinks,breaklinks=true]{hyperref}

\makeatletter
\newcommand\renameautorefname[3][\@empty]{%
  \begingroup
    \def\ult@tmp@langname{\@empty}%
    \ifx#1\@empty
      \global\@namedef{#2autorefname}{#3}%
      \ifcsname languagename\endcsname
        \let\ult@tmp@langname\languagename\relax
      \fi
    \else
      \def\ult@tmp@langname{#1}
      \ifcsname languagename\endcsname
        \ifx#1\languagename
          \global\@namedef{#2autorefname}{#3}%
        \fi
      \fi
    \fi
    \ifx\ult@tmp@langname\@empty\else
      \ifcsname extras\ult@tmp@langname\endcsname
        \expandafter\g@addto@macro\csname extras\ult@tmp@langname\endcsname
          {\@namedef{#2autorefname}{#3}}%
      \fi
    \fi
  \endgroup
}
\makeatother

\renameautorefname{prop}{Proposition}
\renameautorefname{corollary}{Corollary}
\renameautorefname{lemma}{Lemma}
\renameautorefname{definition}{Definition}
\renameautorefname{algorithm}{Algorithm}
\renameautorefname{remark}{Remark}
\renameautorefname{section}{Section}
\renameautorefname{subsection}{Section}

\newcommand{\Y}{\mathcal{Y}} %
\newcommand{\Yn}{\mathcal{Y}_{\textnormal{N}}} %
\newcommand{\Zn}{\mathcal{Z}_{\textnormal{N}}} %
\newcommand{\Yns}{\mathcal{Y}_{\textnormal{\texttt{s}}}} %
\newcommand{\Ynse}{\mathcal{Y}_{\textnormal{\texttt{se}}}}  %
\newcommand{\Ynsne}{\mathcal{Y}_{\textnormal{\texttt{sne}}}}  %
\newcommand{\Ynu}{\mathcal{Y}_{\textnormal{\texttt{u}}}}  %
\newcommand{\Yd}{\mathcal{Y}_{\textnormal{\texttt{d}}}}  %
\renewcommand{\L}{\mathcal{L}} %
\newcommand{\U}{\mathcal{U}} %

\newcommand{\Yl}{\mathcal{Y}_{\lambda}} 
\newcommand{\Ys}{\mathcal{Y}^s} %
\newcommand{\Yg}{\mathcal{G}} %
\newcommand{\mcN}{\textnormal{N}} %
\newcommand{\G}{\mathcal{G}} %
\newcommand{\C}{\mathcal{C}(y)}
\newcommand{\D}{\mathcal{\Y}(y^s)}

\newcommand{\I}{\mathcal{I}} %
\newcommand{\J}{\mathcal{J}} %
\renewcommand{\S}{\mathcal{S}} %

\newcommand{\R}{\mathbb{R}} %
\newcommand{\Rp}{\mathbb{R}^p} %
\newcommand{\Rpg}{\mathbb{R}_{>}^p} %
\newcommand{\Rpgeqq}{\mathbb{R}_{\geqq}^p} %
\newcommand{\Rpgeq}{\mathbb{R}_{\geq}^p} %

\newcommand{\+}{\oplus} %
\newcommand{\MinkSum}{\oplus}

\newcommand{\msum}{\bigoplus} %
\newcommand{\x}{\times} %
\newcommand{\nondom}[1]{ \left(#1\right)_\textnormal{N}}  %
\newcommand{\tp}[1]{#1^{\textnormal{T}}}  %
\newcommand{\hull}[1]{\textnormal{conv}(#1)}  
\renewcommand{\int}[1]{\textnormal{int}(#1)}   %

\newcommand{\nd}{ND\@\xspace} %
\newcommand{\eg}{e.g.\@\xspace} %
\newcommand{\ie}{i.e.\@\xspace} %
\newcommand{\msp}{MSP\@\xspace} %
\newcommand{\NSP}{Decision Space Minkowski Sum Problem\@\xspace} %
\newcommand{\nsp}{DS-MSP\@\xspace} %
\newcommand{\nsfp}{MSP\@\xspace} %
\newcommand{\MSP}{Minkowski Sum Problem\@\xspace} %
\newcommand{\mgs}{MGS\@\xspace} %

\newcommand{\mL}{\texttt{\textsc{l}}\xspace} %
\newcommand{\mU}{\texttt{\textsc{u}}\xspace} %
\newcommand{\mM}{\texttt{\textsc{m}}\xspace} %
\newcommand{\mLU}{\mL{}\mU} %

\begin{document}

\title{Generator Sets for the Minkowski Sum Problem - Theory and Insights\\
}
\cortext[cor1]{Corresponding author}
\author[Aarhus University]{Mark Lyngesen\corref{cor1}}
\ead{lyngesen@econ.au.dk}
\author[Aarhus University]{Sune Lauth Gadegaard}
\ead{sgadegaard@econ.au.dk}
\author[Aarhus University]{Lars Relund Nielsen}
\ead{larsrn@econ.au.dk}
\affiliation[Aarhus University]{organization={Department of Economics and Business Economics, BSS, Aarhus University},%
	    country={Denmark}
}

\journal{arXiv}

\begin{abstract}
This paper considers a class of multi-objective optimization problems known as Minkowski sum problems. Minkowski sum problems have a decomposable structure, where the global nondominated (Pareto) set corresponds to the Minkowski sum of several local nondominated sets.
In some cases, the vectors of local sets does not contribute to the generation of the global nondominated set, and may therefore lead to wasted computational efforts.
Therefore, we investigate theoretical properties of both necessary and redundant vectors, and propose an algorithm based on bounding sets for identifying unnecessary local vectors. We conduct extensive numerical experiments to test the the impact of varying characteristics of the instances on the resulting global nondominated set and the number of redundant vectors. 
\end{abstract}


\maketitle

\section{Introduction}

This paper examines the computation of the nondominated (\nd) set for the Minkowski sum of multiple local sets, denoted \emph{the \nd sum}. Efficient computation of the \nd sum is essential in multi-objective optimization (MO) problems that can be decomposed into subproblems, where each subproblem contributes to the \nd vectors of the problem.
Such problems arise in various contexts, such as organizational decision-making, where departmental decisions impact shared objectives like costs and CO\textsubscript{2} emissions, as well as in complex systems with independent subsystems or coupled subproblems \citep{Gardenghi2011,Kerberenes2022phd}.

Mathematical optimization problems are often easier to solve when decomposed into smaller independent subproblems rather than tackled as a single large problem. In MO, decomposition frequently refers to \emph{objective space decomposition} (\emph{OSD}), where the problem's objective function vector is represented as coordinate-wise combinations of the subproblems objectives. The subproblems may share identical feasible decision space and address a subset of objectives in each subproblem \citep{engau2007} or partition the feasible objective space into subproblems \citep{liu2013}. Furthermore, subproblems can be linked in decision space, as seen in complex and interwoven systems \citep{dietz2020introducing,klamroth2017}. The primary aim of OSD is to simplify the selection process among multiple alternatives or when a complete problem formulation is impractical \citep{dandurand2015}.

One may also consider decomposing the decision space, \ie \emph{decision space decomposition} (\emph{DSD}). DSD focuses on reducing the computational effort for finding the efficient solution set by decomposing the problem. 
We are considering DSD problems where the decision space is decomposed into subproblems, which are connected in objective space. 
Specifically, we consider problems %
where the objective vectors of the problem can be expressed as the vector sum of the subproblems objective vectors. This is referred to as the \emph{\NSP} (\nsp) since the \nd set can be described as the \nd sum of the Minkowski sum of the objective vectors in the subproblems.

\citet{Gardenghi2011} studies algabraic properties of \nsp{s} and show that to solve a \nsp it suffices to compute only the \nd vectors of each subproblem. \citet{Kerberenes2022phd} and \citet{klamroth2024} present a solution approach for \nsp{s} (with binary and general integer variables, respectively) in which the \nd set of each subproblem is computed and thereafter the \nd sum is computed. The last step is referred to as the \emph{\MSP} (\nsfp) since the \nd sum of the \emph{local sets} is found without considering the decision space of each subproblem.

\cite{Kerberenes2022} consider \nsfp{s} with several \nd local sets and propose a sequential and doubling algorithm, which computes the \nd sum by iteratively computing \nd sums of two sets. \cite{hespe2023} consider \nsfp{s} in the bi-objective case with two \nd local sets. By showing that the convex hull of the \nd sum can be computed in linear time using the extreme supported local vectors they present a so-called successive sweep algorithm and apply it to a bi-objective route planning problem. \cite{klamroth2024} present theoretical bounds on the number of elements in the \nd sum of two local sets and report improved running times of the sequential algorithm for \nsfp using an implementation of the divide-and-conquer filtering procedure by \cite{gomes2018}.

Most of the reviewed literature focus on efficient filters for \nsfp{s} as a solution method for solving \nsp{s}. However, often in \nsp{s}, the \nd set of a subproblem is found solving an MO optimization problem, and may therefore be time-consuming to achieve. Moreover, often only a fraction of the local \nd vectors are necessary for generating the \nd sum. If a local \nd vector is not needed to find the \nd sum, then it is \emph{redundant} and computing it can be considered as a waste of computational effort.

In this paper, we consider the \nsfp focusing on finding the \nd sum for the Minkowski sum over local sets containing \nd objective vectors. The paper contributes with several theoretical and computational results which can be summarized as follows:
\begin{enumerate}
    \item We present theoretical results related to the classification of the \nd local vectors and the impact on classification of the vectors in the \nd sum.
    \item We introduce the concept of generator sets needed to compute the \nd sum and study their properties. Here we seek small subsets sufficient for computing the \nd sum or equivalently identifying redundant vectors. 
    \item Two algorithms are given for finding a minimum generator set and a generator set based on bounding sets of each local set, respectively. 
    \item Through computational tests on diverse instances, we analyze how the shape of local sets impacts the \nd sum. Additionally, we demonstrate that, depending on the shape of the local sets, the number of objectives, and the number of subproblems, many vectors become redundant. A simple algorithm based on bounding-sets efficiently identifies many of these redundancies.
\end{enumerate}

The remainder of this paper is structured as follows. \autoref{section:prerequisites} introduces notation and prerequisites for the paper.
Then, \autoref{section:NSP} presents the \nsfp along with theoretical results and algorithms. \autoref{sec:results} presents computational results and conclusions are given in \autoref{sec:conclu}.

\section{Preliminaries}\label{section:prerequisites}

For sets $\Y^1, \Y^2 \subseteq\Rp$, let $\Y^1 \MinkSum \Y^2=\{y\in\Rp \mid y=y^1+y^2,\ y^1\in\Y^1,\ y^2\in \Y^2\}$ denote the \emph{Minkowski Sum} (\emph{MS}) of the two sets. This definition can be extended to multiple sets $\Y^s\subseteq\Rp$, $s\in\mathcal{S} = \{1,...,S\}$, as follows,
\begin{equation}
    \Y = \msum_{s \in \S} \Y^s = \Y^1 \+ \Y^2 \+ \cdots  \+ \Y^S =\{y\in\Rp \mid y=\sum_{s\in S}y^s, y^s\in\Y^s\}
\end{equation}
Note that the Minkowski sum is associative %
and commutative. %

Define the following binary relations to compare vectors $y^1,y^2\in \Rp$:
\begin{align}
    y^1&\leqq y^2 \Leftrightarrow y^1_i\leq y^2_i, &&\forall i\in\{1,..,p\}\\
    y^1&\leq y^2 \Leftrightarrow y^1\leqq y^2, y^1\neq y^2\\
    y^1&< y^2 \Leftrightarrow y^1_i<y^2_i, &&\forall i\in\{1,..,p\}
\end{align}
If $y^1\leq y^2$, we say that $y^1$ \emph{dominates} $y^2$. Similarly, if $y^1\leqq y^2$ ($y^1<y^2$) we say that $y^1$ \emph{weakly dominates} $y^2$ ($y^1$ \emph{strictly dominates} $y^2$). Given $\Y \subseteq \Rp$ let the \emph{\nd set} of $\Y$ be $\Yn = \{y\in\Y \mid \not\exists y'\in\Y : y'\leq y\}$ and the \emph{dominated set} $\Yd = \Y \setminus \Yn$. 
Moreover, a set $\Y$ is \emph{stable} if $\Y=\Yn$.

The binary relations for vectors defined above can naturally be generalized to sets. Let $\Y^1,\Y^2 \subseteq \Rp$ be two stable sets then
\begin{align}
    \Y^1&\leqq \Y^2\quad\Leftrightarrow\quad \forall y^2\in \Y^2\ \exists y^1\in \Y^1: y^1\leqq y^2\ \\
    \Y^1&\leq \Y^2\quad \Leftrightarrow \quad\Y^1\leqq\Y^2, \Y^1\neq \Y^2\\
    \Y^1&<\Y^2\quad \Leftrightarrow \quad \forall y^2\in \Y^2\ \exists y^1\in \Y^1: y^1\leq y^2\ 
\end{align}
Letting $\Rpgeqq = \{y\in\Rp \ | \ y \geqq 0\}$ and define $\Rpgeq$ and $\Rpg$ analogously, we have:
\begin{definition}\label{def:nd}
    Given $\Y \subseteq \Rp$ and $\lambda \in \Rpg$, let $\Yl = \arg\min\{\tp{\lambda}y \mid y\in\Y\}$. %
    The \nd set $\Yn \subseteq \Rp$ can be partitioned into \emph{supported} \nd vectors
	\begin{align}
		\Yns %
             &= \left\{ y \in \Y \mid \exists \lambda : y \in \Yl \right\} \label{eq:yns} = \bigcup_{\lambda\in\mathbb{R}^p_>}\mathcal{Y}_\lambda,
	\end{align}
    and \emph{unsupported}, $\Ynu = \Yn \setminus \Yns$, \nd vectors.
    Furthermore, $\Yns$ can be partitioned into (supported) \emph{extreme} 
	\begin{align}
		\Ynse &= \left\{ y \in \Y \mid \exists \lambda \in \Rpg : y \in \Yl, |\Yl| = 1 \right\}, \label{eq:ynse1}
	\end{align}
    and \emph{supported non-extreme} $\Ynsne = \Yns \setminus \Ynse$. 
        
\end{definition}
An illustration of the different classifications of $\Y$ is given in \autoref{fig:msp-ex}.

\section{The Minkowski Sum Problem}\label{section:NSP}

In the following, let $\mathcal{S}=\{1,...,S\}$ be an index set, and consider finite non-empty local sets $\Ys \subseteq \R^p$.
The \emph{\MSP} (\emph{MSP}) can be stated as 
\begin{equation}
\min\{y \mid y \in \Y = \msum_{s \in \S} \Ys  \}.\label{eq:MinkSumProb}
\end{equation}
Note, the MSP is multi-objective in nature, since solving \eqref{eq:MinkSumProb} should be interpreted as finding the \nd set $\Yn$, of the Minkowski sum $\Y$. Note that $\Yn$ is unique and denoted the \emph{\nd sum} since each vector $y = y^1 + \ldots + y^S, y^s\in\Ys, y\in\Y$ is the \emph{vector sum} (\emph{VS}) of local vectors.

\begin{figure}[tb]
    \centering
    \resizebox{0.9\linewidth}{!}{\input{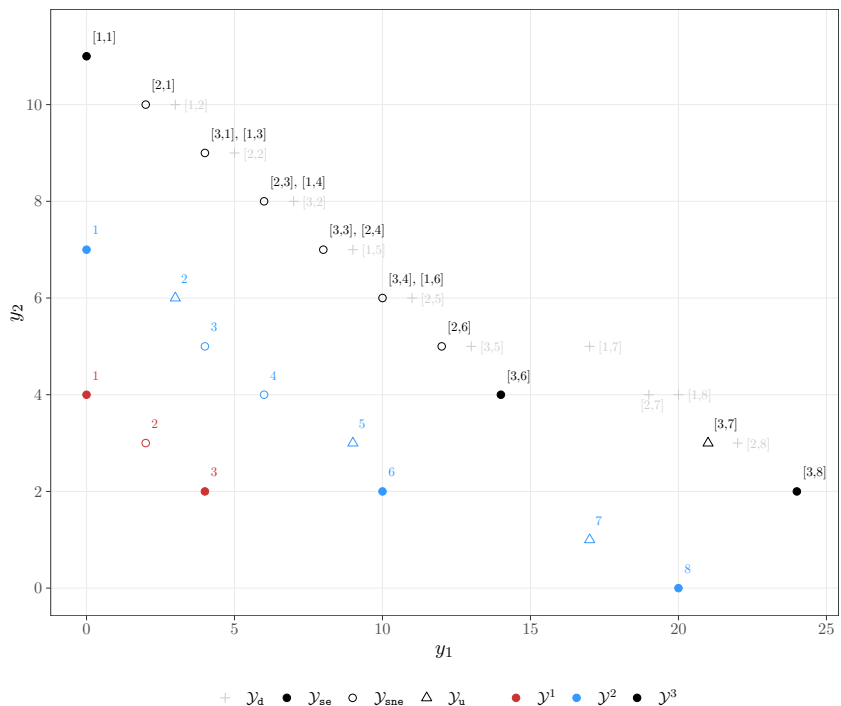}}
    \caption{Classifications of two local sets and their MS. The number/index besides each local vector is used to identify which combinations are used to generate a VS, \eg VS $(10, 6)$ is the sum of the vector with index 3 from $\Y^1$ and the vector with index 4 from $\Y^2$. Moreover, index $[1,6]$ gives the same VS. Vectors in grey are dominated.} \label{fig:msp-ex}
\end{figure}

In general, not all local vectors are needed to find $\Yn$. It is well known, that in order to find $\Yn$, only the \nd set of each local set is necessary \citep{Gardenghi2011}. That is $\Yn \subseteq \msum_{s \in \S} \Yn^s \label{eq:msyn}$.
However, often not all vectors in $\Yn^s, s\in\S$ are needed either. An example is given in \autoref{fig:msp-ex}. Here, the VSs obtained using the vector $y^2\in\Yn^2$ are dominated. That is, $y^2$ is not needed for finding $\Yn$.

\subsection{Relations between \nd vectors}

Some important relations exists depending on the classification of the \nd vector. First, the focus will be on the supported \nd sets.

\begin{prop}\label{prop:s}
    Consider local sets $\Ys$, for $s\in\S$, and let $\Y = \msum_{s \in \S} \Y^s$. Given $\lambda\in\Rpg$, let $\Yl = \arg\min\{\tp{\lambda}y \mid y\in\Y\}$ and $\Yl^s = \arg\min\{\tp{\lambda}y \mid y\in\Y^s\}$.  Then the following relations hold
    \begin{enumerate}
        \item If $y = y^1 + \ldots + y^S \in\Yns $ then $ y^s\in\Yns^s,  $ for all $s\in\S$.\label{sup1}
        \item $\Yl = \msum_{s \in \S} \Yl^s = (\Yl)_{\textnormal{N}} \subseteq \Yns$.\label{sup2}
    \end{enumerate}
\end{prop}

\begin{proof} \ %
\begin{enumerate}

\item If $y\in\Yns$, then due to \eqref{eq:yns} we can find $\lambda \in \Rpg$ satisfying 
\begin{equation}
    \begin{aligned}
        \tp{\lambda} y &= \min\{\tp{\lambda} \bar{y} \mid \bar{y}\in\Y\} \\
            &= \min\{\tp{\lambda}(\bar{y}^1 + \ldots + \bar{y}^S) \mid \bar{y}^s\in\Ys, s\in \S\} \\
            &= \sum_{s\in\S} \min\{\tp{\lambda}\bar{y}^s \mid \bar{y}^s\in\Ys\}  
    \end{aligned}\label{eq:lmin}
\end{equation}
Now, assume for contradiction that there exists an $\hat{s}\in\S$ such that $y^{\hat{s}}\not\in\Yl^{\hat{s}}$. This, in turn, implies that there exists a $\tilde{y}^{\hat{s}}\in\Yl^{\hat{s}}$ such that $\lambda^T\tilde{y}^{\hat{s}}<\lambda^T{y}^{\hat{s}}$, contradicting \eqref{eq:lmin}. Hence, $y^{\hat{s}}\in\Yl^{\hat{s}}$, and thereby $y^{\hat{s}}\in \Yns^{\hat{s}}$.

\item From the proof of 1. we have that $y = y^1 + \ldots + y^S \in \Yl$ implies $y^s\in \Yl^s$. Furthermore, using \eqref{eq:lmin} we see that $y^s\in \Yl^s$ implies  $y\in \Yl$. Hence, $\Yl=\msum_{s\in\S}\Yl^s$. That $\Yl =(\Yl)_{\textnormal{N}}$ and $\Yl \subseteq \Yns$ follows from \autoref{def:nd}.
\end{enumerate}
\end{proof}

\autoref{prop:s} states that supported vectors in $\Yn$ are obtained as a VS of supported local vectors $\Yn^s$ (Relation \ref{sup1}). These can be found by solving scalarized problems for each local set (Relation \ref{sup2}). That is, by selecting a \emph{search direction} $\lambda$. Note that the search direction $\lambda$ in Relation \ref{sup2} must be fixed. The vector sum of supported local vectors for different search directions does not necessarily give a supported vector, as can be seen in \autoref{fig:msp-ex}. Here the VS of $y^1\in\Yns^1$ and $y^8\in\Yns^2$ is dominated. Stronger relations can be obtained by considering the extreme \nd set.

\begin{prop}\label{prop:se}
    Consider $ y\in\Y$ where $ y = y^1 + \ldots + y^S$ with $y^s\in\Ys$ for $ s\in\S$. Then the following relations hold
    \begin{enumerate}
        \item If $\lambda\in\Rpg $ and $\Yl = \{ y \}  $ then $ y\in\Ynse$ and $y^s\in \Ynse^s$ for all $s\in\S$.\label{sup3}
        \item If $y \in\Ynse $ then $ y^s\in\Ynse^s$ for all $s\in\S$.\label{sup4}
        \item If $y \in\Ynse$ then no other VS equals $y$ (uniqueness).\label{sup5}
        \item Let $y^{\hat{s}}\in\Ynse^{\hat{s}}$ for some $ \hat{s}\in\S$. Then there exists a extreme vector, $y\in\Ynse$, such that $y=y^1+\dots+y^{\hat{s}}+\dots+y^S$, where $y^s\in\Ynse^s$ for all $s\in\S$.\label{sup6}
    \end{enumerate}
\end{prop}

\begin{proof} 
Let us first prove that for $\lambda \in \Rpg$
\begin{equation}
    |\Yl| = 1 \Leftrightarrow |\Yl^s| = 1, \forall s\in\S \label{eq:se-l}
\end{equation}
$\Rightarrow$: Assume for contradiction that (wlog) local set $\Yl^1$ satisfies that $\hat{y}^1, \bar{y}^1  \in \Yl^1, \hat{y}^1\neq\bar{y}^1$, then picking the same $y^s\in\Yl^s, s = 2, \ldots, S$ and letting $\hat y = \hat y^1 + \sum_{s = 2}^S y^s$ and  $\bar y = \bar y^1 + \sum_{s = 2}^S y^s$ we have 
\begin{align}
    \tp{\lambda}\hat y = \tp{\lambda}\hat y^1 + \sum_{s = 2}^S y^s = \tp{\lambda} \bar y^1 + \sum_{s = 2}^S y^s = \tp{\lambda}\bar y %
\end{align}
which contradicts that $|\Yl| = 1$ since $\hat{y}, \bar y \in \Yl$ and $\hat y \neq \bar y$. 

\noindent $\Leftarrow$: Assume for contradiction that $|\Yl| > 1$. Then $\exists\Yl^s: |\Yl^s|>1$ which is a contradiction.%

\begin{enumerate}
\item If $\Yl = \{y\} $ then $ y \in\Ynse$ due to \eqref{eq:ynse1}. Moreover, due to \eqref{eq:se-l} $|\Yl^s| = 1 \Rightarrow y^s \in\Ynse^s, \forall s \in \S$. 

\item This is a consequence of using \eqref{eq:ynse1} and Relation \ref{sup3}.

\item To prove uniqueness, pick $\lambda$ such that $y \in \Yl$ satisfies \eqref{eq:se-l} and assume there exists a $\hat{\lambda} \neq \lambda $ such that $ y = \hat{y}^1 + \ldots + \hat{y}^S \in \Y_{\hat{\lambda}}$ and (wlog) $\hat{y}^1 \neq y^1$. Then $\tp{\lambda} y^1 < \tp{\lambda} \hat{y}^1$ and $\tp{\lambda} y^s \leq \tp{\lambda} \hat{y}^s$ for $ s=2,\ldots,S$ and
\begin{align}
    \tp{\lambda} y &= \tp{\lambda} (y^1 + y^2 + \ldots + y^S) \\
              &< \tp{\lambda} (\hat{y}^1 + y^2 + \ldots + y^S) \\
              &\leq \tp{\lambda} (\hat{y}^1 + \hat{y}^2 + \ldots + \hat{y}^S) = \tp{\lambda} y,
\end{align}
which is a contradiction.

\item If $y^{\hat{s}}\in\Ynse^{\hat{s}} $ then there exists a  $\hat{\lambda}$ with $ \Y_{\hat{\lambda}}^{\hat{s}} = \{y^{\hat{s}}\}$. Due to Relation~\ref{sup2} in \autoref{prop:s} all $y\in\Ynse \cap \Y_{\hat{\lambda}}$ must use $y^{\hat{s}}$ in their VS.
\end{enumerate}
\end{proof}

\autoref{prop:se} states that if a vector is extreme in the \nd sum, then it is obtained as the sum of extreme supported local vectors (Relation \ref{sup3}-\ref{sup4}).
Moreover, an extreme VS is obtained as a unique sum of local vectors (Relation \ref{sup5}). In addition, all extreme local vectors are needed for finding the extreme vectors in the \nd sum (Relation~\ref{sup6}). 

Note \autoref{prop:se} does not state that the VS of extreme vectors give an extreme vector as can be seen in \autoref{fig:msp-ex}. Here, the VS of extreme local vectors $y^1\in\Yns^1$ and $y^8\in\Yns^2$ is dominated. This even does not hold if $\Yl$ contains multiple solutions, \eg by setting $\lambda = (0.5, 0.5)$ we get $y^3\in\Yl^1$ and $y^1\in\Yl^2$. Both are extreme, but the VS gives a supported non-extreme vector. 

The largest set of extreme local vectors provide a lower bound on the number of extreme vectors in the \nd sum.
 This follows from Relation \ref{sup6} as each $y^{\hat s} \in \Ynse^{\hat s}$ is part of the VS of some $y \in \Ynse$. Also, each such $y \in \Ynse$ has a unique VS (Relation \ref{sup5}). This gives us the following corrollary:
\begin{corollary} 
        $|\Yn| \ge |\Ynse| \ge \max_{s \in \S}|\Ynse^s|$. \label{sup7}
\end{corollary}
It does not generally hold that $|\Yn| \ge \max_{s\in \S}|\Yn^s|$ as shown in \cite{stiglmayr2014}. Now consider the supported non-extreme set.

\begin{prop}\label{prop:sne}
    Consider $y=y^1 + \ldots + y^S \subseteq \Y = \msum_{s \in \S}\Y^s$ and assume that $y^s \in \Yns^s$. If $y^{\hat s} \in \Ynsne^{\hat s}$ for some $\hat s \in \S$ then 
    $y \in \Ynsne$.
\end{prop}
\begin{proof} 
Note if $y^{\hat{s}}\in \Ynsne^{\hat{s}}$ it can be written as a convex combination of vectors $\{ \hat{y}^1, \ldots, \hat{y}^k \} \subset \Yl^{\hat{s}}$ for some $\lambda \in \Rp$, \ie
\begin{equation}
    y^{\hat{s}} = \sum_{i = 1}^{k} \alpha_i \hat{y}^i,\quad \alpha \in \R^{k}_\ge, \tp{1}\alpha = 1. \label{eq:ynsne}
\end{equation}
This implies that 
\begin{align}
    y &= y^{1} + \ldots + \sum_{i = 1}^{k} \alpha_i \hat{y}^i + \ldots + y^{S} \\
      &= \sum_{i = 1}^{k} \alpha_i y^1 + \ldots + \sum_{i = 1}^{k} \alpha_i \hat{y}^i + \ldots + \sum_{i = 1}^{k} \alpha_i y^S \\
      &= \sum_{i = 1}^{k} \alpha_i (y^1 + \ldots + \hat{y}^i + \ldots + y^S) = \sum_{i = 1}^{k} \alpha_i \bar{y}^i 
\end{align}
Since $\hat{y}^i\in\Yl^{\hat{s}} \Rightarrow \bar{y}^i = y^1 + \ldots + \hat{y}^i + \ldots + y^S\in\Yns$ and then $y\in\Ynsne$ due to \eqref{eq:ynsne}.
\end{proof}

\autoref{prop:sne} states that a VS of supported \nd vectors where one vector is non-extreme cannot be extreme.

\subsection{Generator sets}\label{subsection:generator-sets}
Although only $\Yn^s, s\in\S$ is needed to find the \nd sum, it may happen that not all vectors in $\Yn^s$ are necessary. In particular, there may exist local vectors that 
always lead to dominated vectors. For instance, consider the example given in \autoref{fig:msp-ex}. 
Here, any vector in the MS of $\Yn^1$ and $\{y^2, y^5\} \subseteq \Yn^2$ is dominated.
Such vectors are undesirable for two obvious reasons: first, in case the local sets, $\Yn^s$, are known these vectors add unnecessary computational work when computing the \nd sum. Second, if the \nd set of the local sets are not known a priori, but instead emerge from a multi-objective optimization problem, then computational effort is wasted when computing such vectors. For this reason, we introduce the concept of generator sets
\begin{definition}
    Given stable sets $\G^s\subseteq \Y^s$ with $\Yn=\left(\msum_{s\in \S} \G^s \right)_\textnormal{N}$ we say that $\G^s$ is a \emph{generator} and that $\G=\{\G^1,...,\G^S\}$ is a \emph{generator set}
\end{definition}
Multiple generator sets may exist and trivially $\G = \left\{\Yn^1,\ldots, \Yn^S \right\}$ is a generator set. Naturally, we seek small generator sets.
\begin{definition}
    A generator set $\G$ is a \emph{minimal generator set} if there does not exist $\tilde{\G} = \{ \tilde\Yg^1,...,\tilde\Yg^S\}$ with $\tilde{\G}^s\subseteq {\G}^s, s\in\S$, with at least one strict inclusion, such that  $\Yn=\left(\msum_{s\in \S} \G^s \right)_\textnormal{N}$ is satisfied.
\end{definition}

Given a generator set $\G$, a vector $y^s\in\G^s$ is \emph{redundant} if $\hat\G = \{ \Yg^1,\ldots,\G^s\setminus\{y^s\},\ldots,\Yg^S\}$ is a generator set. That is, a minimal generator set contains no redundant vectors.
Note, a minimal generator set is not uniquely defined. For instance, consider the example given in \autoref{fig:msp-ex}. Here, all vectors in $\Yn^1$ are needed to find the \nd sum, \ie $\Yg^1 = \Yn^1$. However, for generator $\Yg^2$ both $\{ y^1, y^3, y^6, y^7, y^8 \}$ and $\{ y^1, y^4, y^6, y^7, y^8 \}$ can be used to define a minimal generator set.

If the generator set $\G$ is found by solving a multi-objective optimization problem, it is of interest to keep the set of computed solutions as small as possible. For that reason, we also introduce the concept of minimum (cardinality) generator sets

\begin{definition}\label{def:MGS}
    We say that a generator set $\G = \{ \Yg^1,...,\Yg^S\}$ is a \emph{minimum generator set} if $\sum_{s\in\S} \vert\Yg^s\vert \leq \sum_{s\in\S} \vert \tilde{\Yg}^s\vert$ for all generator sets $\tilde\G = \{ \tilde\Yg^1,...,\tilde\Yg^S\}$.
\end{definition}

A minimum generator set is, in principle, a generator set where the MS is the least computationally demanding to generate. The following lemma states the relation between minimal and minimum generator sets:
\begin{prop}\label{prop:minimumIsMinimal}
    A minimum generator set is a minimal generator set.
\end{prop}
\begin{proof}
    If the generator set $\G=\{\G^1,\dots,\G^S\}$ is not minimal, then (wlog) there exist a proper subset $\tilde{\G}^1\subset \G^1$ such that $\G=\{\tilde{\G}^1,\G^2,\dots,\G^S\}$ is a generator set, implying $\G$ is not a minimum generator set.
\end{proof}
Note, a minimal generator set is not necessarily a minimum generator set.
To identify redundant vectors we introduce the following notation:
\begin{align}
    \C &= \{ (y^1,\ldots, y^S) \mid y = y^1 + \ldots + y^S, y^s\in\Ys, s\in\S \} \label{def:combinationsC}, \forall y \in \Y\\
    \D &= \{y^s\} \+ \left(\msum_{\hat{s}\in\S\setminus\{s\}} \Y^{\hat{s}} \right), \forall y^s \in \Y^s, \forall s \in \S %
\end{align}
The set $\C$ contains all combinations of local vectors with VS equal $y$ and the set $\D$ is the MS obtained given a fixed vector $y^s$ from local set $s$. 

\begin{prop}\label{prop:gen-clas}
    Consider $y = y^1 +\ldots + y^S\in\Y$ with $y^s \in \Y^s, \forall s \in \S$. Then 

    \begin{enumerate}
        \item If $y\in\Yn$ and $|\C|=1$ then any generator set $\G$ have $y^s \in \G^s, \forall s\in\S$. \label{card1}
        \item If $y\in\Yn$ and $y^{\hat s} = c^{\hat s}$ (fixed to $c^{\hat s}$) for all $(y^1, \ldots, y^{\hat s}, \ldots, y^S)\in\C$ then any generator set $\G$ have $y^{\hat s} \in \G^{\hat s}$.\label{card-more}
        
        \item If $\D\subseteq\Yd$ then $y^s\notin\G^s$ for any minimal generator set $\G$. \label{card-dom}
    \end{enumerate}
\end{prop}

\begin{proof}
    First, note that Relation \ref{card1} is a consequence of Relation \ref{card-more}. 

    \begin{enumerate}[start=2]
        \item Assume wlog that $\exists\G : y^1 \notin \G^1$ but then $\exists \hat{y}^1 : y = \hat{y}^1 + y^2 + \ldots + y^S$ implying that $(\hat{y}^1, y^2, \ldots, y^S)\in\C$ which contradicts $y^1 = c^1, \forall c \in \C$.

        \item Assume that there exists a minimal generator set $\G$ with $y^s\in\G^s$, but then $\exists y\in\D : y\in\Yn$ which contradicts $\D\subseteq\Yd$.
    \end{enumerate}
\end{proof}

\autoref{prop:gen-clas} states that if an \nd VS is obtained using unique local vectors, then all these vectors cannot be removed from a generator set (Relation \ref{card1}). This also holds even if more combinations gives the same VS. In that case, if a local vector is used in all combinations, then it must be part of any generator for the corresponding local set (Relation \ref{card-more}). Finally, if all combinations given a local vector are dominated, then it is redundant (Relation \ref{card-dom}).

The link between generator sets and the extreme vectors are readily established by combining Relation \ref{sup5} in \autoref{prop:se} and Relation \ref{card1} in \autoref{prop:gen-clas} as a corollary.
\begin{corollary}\label{cor:extremeAreGenerators}
    For any generator set $\G=\{\G^1,\dots,\G^S\}$, $\Ynse^s\subseteq \G^s$ for all $s\in\S$.
\end{corollary}

An implication of \autoref{cor:extremeAreGenerators} is that if $\Yn^s$ only contains extreme vectors, then $\Yg^S$ is uniquely defined by $\Ynse^s$ in any generator set.

\subsection{Identifying minimum generator sets}

Using \autoref{def:MGS}, a minimum generator set is an optimal solution $\G^* = ({\G^{1}}^*, \ldots, {\G^{S}}^*)$ to the optimization problem
\begin{align}
\begin{split}
     \min_{\G}& \sum_{s \in \S} \lvert \G^s \rvert  \\
	s.t. \quad & \Yn \subseteq \msum_{s \in \S}  \G^s,  \\
		   & \G^s \subseteq \Y^s, \quad \forall s \in \S.
\end{split}
\label{prob:MGS}
\end{align}
    
This problem can be solved using a combinatorial programming formulation with binary variables related to whether $y^s\in\G^s$ or not. However, the solution time can be improved by initially finding local vectors $y^s\in\Ys$ that must be in ${\G^{s}}^*$.

Relation \ref{card-more} in \autoref{prop:gen-clas} states that if a vector $y \in \Yn$ use a local vector $y^s \in \Y^s$ for any VS, then $y^s \in {\G^{s}}^*$. That is, 
\begin{equation}
    \bar{\Y}^s := \{ y^s \in \Y^s \vert \exists y \in \Yn : y^s = c^s, \forall (c^1, \ldots, c^s, \ldots, c^S)\in \C \} \subseteq {\G^{s}}^*.
\end{equation}
This implies that $\Ynse^s \subseteq \bar{\Y}^s$ for all $s \in \S$ due to \autoref{cor:extremeAreGenerators}.

Furthermore, following Relation 3 in \autoref{prop:gen-clas} only local vectors that are part of at least one VS of vectors in the \nd sum need to be considered. That is, ${\G^{s}}^* \subseteq \hat{\Y}^s := \{ y^s \in \Y^s \vert  (\D)_\text{N}\cap \Yn \neq \emptyset \} $.

Hence $\bar{\Y}^s \subseteq {\G^{s}}^* \subseteq \hat{\Y}^s$ and since $\bar{\Y}^s$ is unique, we have the following corollary.

\begin{corollary}
    \label{col:reducesEqualFixed}
    Consider an optimal solution $\G^* = ({\G^{1}}^*, \ldots, {\G^{|S|}}^*)$ to \eqref{prob:MGS} then 
    \begin{enumerate}
        \item If $\bar{\Y^s} = \hat{\Y^s}$ then ${\G^{s}}^* = \bar{\Y^s}$ is a unique generator for local set $s$.
        \item If $\bar{\Y^s} = \hat{\Y^s}$ for all $s \in \S$ then $\G^*$ is a unique minimal and minimum generator set (a unique solution).
        \item If $\Yn = (\msum_{s \in \S} \bar{\Y^s})_\textnormal{N}$ then $\G^*$ is a unique minimal and minimum generator set (a unique solution).
    \end{enumerate}
\end{corollary}

\begin{algorithm}[tb]
    \caption{Finding a minimum generator set.}
    \KwData{$\left(\Y^1, ..., \Y^S\right)$}
    \KwResult{a mimimum generator set $\G^* = (\G^{1*}, ..., \G^{S*})$ }
    \label{alg:MGS} 
    \BlankLine
    Compute $\Yn$ using a filtering algorithm and store $\C, \forall y\in\Yn$\; \label{l:start}
    Compute the sets $\bar{\Y^s}$ and $\hat{\Y^s}, \forall s \in \S$\;
    \lIf{$\hat{\Y^s} = \bar{\Y^s}, \forall s \in \S$}{\label{alg:MGS-if-fixed-is-reduced}\Return $\G^* = (\bar{\Y^1},\ldots, \bar{\Y^S})$}
    \lIf{$\Yn = (\msum_{s \in \S} \bar{\Y^s})_\textnormal{N}$}{\Return $\G^* = (\bar{\Y^1},\ldots, \bar{\Y^S})$ \label{l:end}}
    \lElse{\Return an optimal solution $\G^*$ to \eqref{prob:MGS}\label{l:call-solver}}
\end{algorithm}

\autoref{alg:MGS} present pseudocode for finding a minimum generator set. First, the sets used to check \autoref{col:reducesEqualFixed} are calculated and used (lines \ref{l:start}-\ref{l:end}). Finally, if \autoref{col:reducesEqualFixed} cannot be applied, a combinatorial optimization problem of \eqref{prob:MGS} is solved. An integer programming formulation is given in \ref{app:IP-MGS}.

\subsection{Identifying a reduced cardinality generator set}

This section considers an algorithm for removing local vectors that cannot be part of any minimal generator set. Assume that each local set $s$ is bounded by two stable sets $\L^s$ and $\U^s$ such that $\L^s \leqq \Yn^s \leqq \U^s$. We refer to the sets $\L^s$ and $\U^s$ as \emph{bounding sets} for $\Y^s$ similar to the notion of bound sets introduced in \citet{ehrgott2007bound}, but is instead based on the binary set relations defined in \autoref{section:prerequisites}. For instance, $\L^s$ may have been found using an outer approximation algorithm and $\U^s$ as a subset of $\Y^s$.

\begin{prop}\label{prop:MSP-lb-ub}
	Consider a MSP with \nd sum $\Yn$ and local sets $\Y^s$, $s\in\S$. Assume that $\L^s \leqq \Yn^s \leqq \U^s$ for all $s \in \S$. Then
	\begin{align}
		\nondom{\msum_{s \in \S} \L^s} \leqq
		\Yn \leqq
		\nondom{\msum_{s \in \S} \U^s}
	\end{align}
\end{prop}
\begin{proof}
    We start by showing the first inequality
		$\nondom{\msum_{s \in \S} \L^s} \leqq \Yn $. Since $\L^s \leqq \Yn^s$ we have
	$\forall y^s \in \Yn^s, \exists l^s \in \L^s: l^s \leqq y^s$ which imples 
	$\forall y = y^1 + \ldots + y^S \in \Yn, \exists l = l^1+ \ldots + l^S \in \msum_{s \in \S} \L^s: l \leqq y \implies \nondom{\msum_{s \in \S} \L^s} \leqq \Yn$.

	Likewise 
	$\forall u^s \in  \U^s, \exists y^s \in \Yn^s: y \leqq u  $
	implies 
	$\forall u = u^1 + \ldots + u^S \in \msum_{s \in \S} \U^s, \exists y = y^1 + \ldots + y^S \in \Yn: y \leqq u $
	showing the second inequality 
	$\Yn \leqq \nondom{\msum_{s \in \S} \U^s}$.
\end{proof}

\autoref{prop:MSP-lb-ub} states that local bound sets can be used to calculate bound sets for the \nd sum. 
Since $\Yn^s \leqq \Yn^s$, we may use the \nd set as bounds for some local sets in \autoref{prop:MSP-lb-ub}.

\begin{corollary}\label{coll:MSP-lb-ub}
	Assume that $\L^s \leqq \Yn^s \leqq \U^s$ for all $s \in \S$. Consider two partitions $\S^1 \cup \S^2 = \S$ and $\S^3 \cup \S^4 = \S$ satisfying $\S^1 \cap \S^2 = \S^3 \cap \S^4 = \emptyset$. Then
	\begin{align}
	    \nondom{\left(\msum_{s \in \S^1} \L^s\right) \+ \left(\msum_{s \in \S^2} \Yn^s\right) }
	    \leqq \Yn \leqq
	    \nondom{\left(\msum_{s \in \S^3} \Yn^s\right) \+ \left(\msum_{s \in \S^4} \U^s\right) }
	\end{align}
\end{corollary}

\begin{prop}\label{prop:conditionalDom}
Assume $\L^s \leqq \Yn^s \leqq \U^s$. 
If for some $\bar s \in \S$ there exists a vector $y^{\bar s} \in \Y^{\bar s}$ which satisfies the set inequality
$$ \nondom{\msum_{s \in \S} \U^s} < \nondom{\{y^{\bar s}\} \+ \left(\msum_{s \in \S \setminus \{\bar s\}} \L^s\right)},$$
then $y^{\bar s} \notin \G^{\bar s}$ for any minimal generator set $\G$.
\end{prop}

\begin{proof}
    Let $\U = \nondom{\msum_{s \in \S} \U^s}$ and $\hat{\L} = \nondom{\{y^{\bar s}\} \+ \left(\msum_{s \in \S \setminus \{\bar s\}} \L^s\right)}$. Since $\hat \L \leqq \nondom{\Y(y^{\bar s})}$ and $\Yn \leqq \U$ (\autoref{prop:MSP-lb-ub}) we have $\Yn \leqq \U < \hat{\L} \leqq \nondom{\Y(y^{\bar s})}$ and hence $\Y(y^{\bar s}) \subseteq \Yd$, \ie $y^{\bar s} \notin \G^{\bar s}$ for any minimal generator set $\G$ (\autoref{prop:gen-clas}, Relation \ref{card-dom}).
\end{proof}

\autoref{prop:MSP-lb-ub} and \autoref{prop:conditionalDom} can be used recursively on local sets. In other words, they apply to any local set $\hat{\S} \subseteq \S$. Specifically, they are valid for all pairs of local sets. 

\begin{algorithm}[tb]
    \caption{Finding reduced cardinality generators using bounding sets}\label{alg:Alg3Pair}
    \KwData{$\L^s \leqq \Yn^s \leqq \U^s, \forall s\in\S$.}
    \KwResult{A generator set $\hat \G = \left(\hat\G^1, \ldots ,\hat\G^S \right)$.}
    \BlankLine
    \ForAll{$\bar s \in \S$\label{alg:Alg3Pair-choose-s} }{

	    let $\hat \G^{\bar s} = \Yn^{\bar s}$ \;
	    \ForAll{$s \in \bar\S = \S \setminus \{\bar s\}$\label{alg:Alg3Pair-S-bar}}{
		    \ForAll{$y^{\bar s} \in \hat\G^{\bar s}$}{
			    \If{$\nondom{\U^{\bar s} \+ \U^{s}} < \nondom{\{y^{\bar s}\} \+ \L^{s}}$ \label{alg:Alg3Pair-if-dominated} }{
				    let $\hat \G^{\bar s} := \hat \G^{\bar s} \setminus \{y^{\bar s}\}\; $
			    }
		    }
	    }
    }
    \Return{$\hat\G = (\hat \G^1,\ldots, \hat \G^S)$}
\end{algorithm}

An algorithm capable of removing local vectors from the current generator set, which cannot be part of any minimal generator set, is described in \autoref{alg:Alg3Pair}. The algorithm employs bounds $\L^s \leqq \Yn^s \leqq \U^s$ for every local set. For each local set $\bar s$, the generator is initialized by assigning $\hat\G^{\bar s} = \Yn^{\bar s}$, and pairwise comparisons are then conducted on each local set $s \in \bar\S \coloneqq \S \setminus \{\bar s\}$ using vectors $y^{\bar s} \in \G^{\bar s}$. If a vector $y^{\bar s}$ meets the criterion of \autoref{prop:conditionalDom}, it is excluded.

\section{Computational results}\label{sec:results}

The purpose of this computational study is to answer the following questions: 

\newcommand{\rqF}{What is the cardinality of the \nd sums and how does this depend on the shape of the \nd local sets?\xspace}
\newcommand{\rqS}{How large is the minimum generator sets?\xspace}
\newcommand{\rqT}{How many redundant vectors can be identified using bounding sets?\xspace}
\begin{enumerate}
    \item \rqF
    \item \rqS
    \item \rqT
\end{enumerate}

All computational results were obtained using a Linux (CentOS 7) cluster with 3.1 GHz (AMD EPYC 9554) cores. To determine the \nd sum $\Yn$ for each \msp instance, we used a C implementation of the LimMem nondominance filter from \cite{klamroth2024} with a memory limit of 64 GB. The LimMem algorithm uses a buffer to iteratively build the \nd sum using subsets of the local sets thereby avoiding storing the entire Minkowski sum in memory. The remaining algorithms were implemented in Python 3.11 and the code can be accessed at \cite{MSP-Generators-Lyngesen24}. A time limit of 12 hours was used for each instance when running the algorithms.

\subsection{Test instances}

\begin{figure}[tbp]
    \centering
    \includegraphics[width=0.99\linewidth]{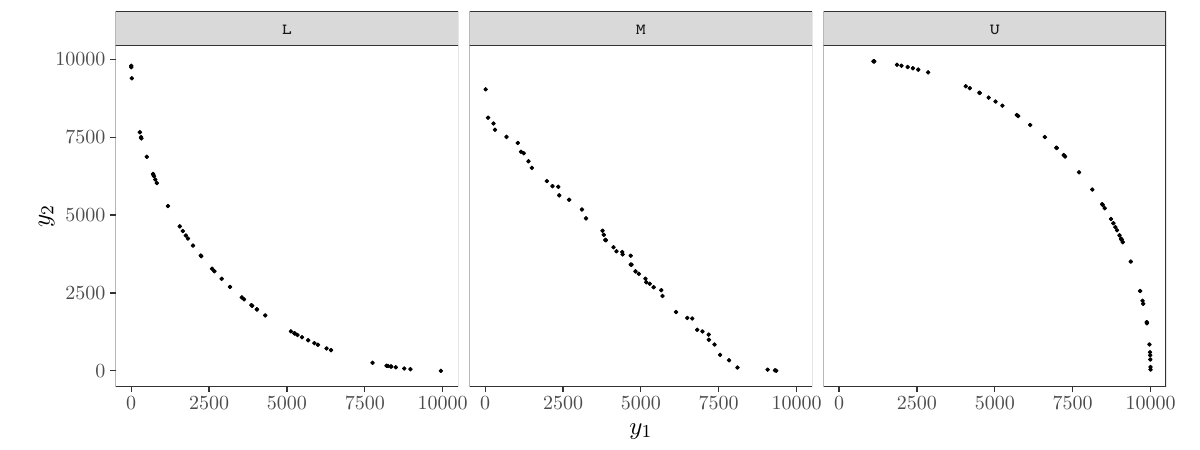}
    \caption{local sets generated using the three different methods \mL, \mM and \mU ($p=2$).} \label{fig:sp-50-lmu}
\end{figure}    

\begin{figure}[tbp]
    \includegraphics[width=0.99\linewidth]{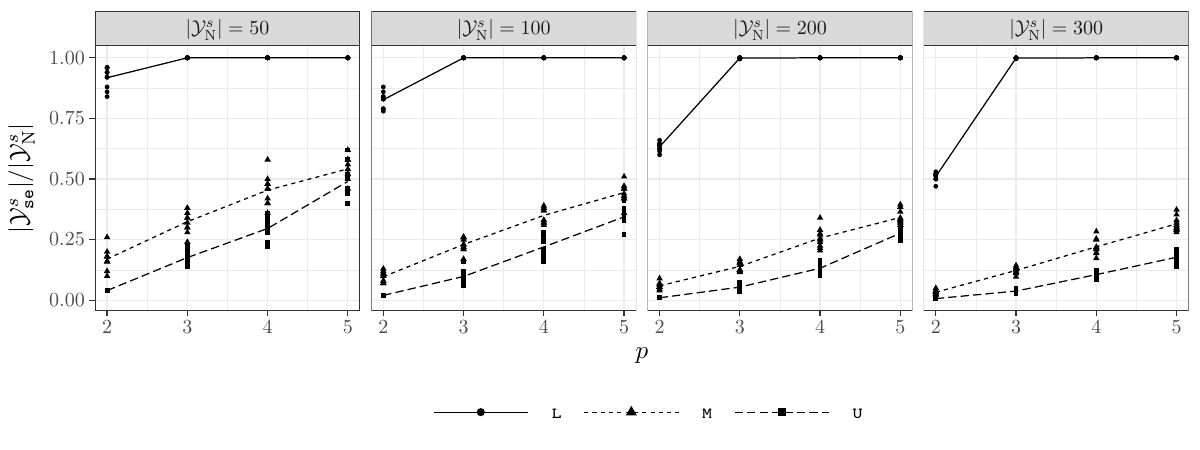}
\caption{Relative number of extreme vectors given number of objectives under different cardinalities of the local sets.}\label{fig:se}
\end{figure}

For the analysis, a collection of local set instances is considered, which are combined into MSP instances. For each local set instance, the \nd set $\Yn^s\in\mathbb{Z}^p$ is generated. The goal is to generate diverse local sets with respect to cardinality ($|\Yn^s|$), number of objectives ($p$) and shape (relative number of extreme vectors). The \nd vectors are generated inside a hypercube $[0, 10000]^p$ until the specified cardinality is obtained. Note that local sets $\Y^s$ can be shifted with a fixed vector $z^s$: define $\Zn^s = \Yn^s \MinkSum \{z^s\}$ and observe that
\begin{align}
    \Zn &= \left( \msum_{s\in\S} \left( \Yn^s \MinkSum \{z^s\} \right) \right)_\mcN \\
        &= \left( \msum_{s\in\S} \Yn^s \MinkSum \{z^1 + \ldots + z^S\} \right)_\mcN \\
        &= \left( \msum_{s\in\S} \Yn^s \right)_\mcN \MinkSum \{ z \} \\
        &= \Yn \MinkSum \{ z \},
\end{align}
where $z = z^1 + \ldots + z^S$. Hence, testing instances generated in hypercubes $[l,u]^p$ with $l\neq 0$ is irrelevant for the computational study.

Three \emph{methods} affecting the shape of the \nd local sets are used. The method \mL generates \nd vectors in the southwest part of a sphere, resulting in many extreme vectors (see \mL in \autoref{fig:sp-50-lmu}). Similarly, the method \mU generates vectors on the northeast part of a sphere, resulting in many unsupported vectors (see \mU in \autoref{fig:sp-50-lmu}). %
Finally, method \mM generates vectors between two parallel hyperplanes with a normal vector pointing towards the origin, resulting in a situation with unsupported vectors near the supported vectors which lies close to one of the hyperplanes (see \mM in \autoref{fig:sp-50-lmu}). Note that all generated vectors are rounded to integers, some vectors may be unsupported even for method \mL. %
Ten instances were generated for each $p=2,\ldots, 5$, cardinality $|\Yn^s| = 50, 100, 200, 300$ and method \mU, \mL and \mM, resulting in 480 local sets in total. For further details see \cite{MOrepo-Lyngesen24}.

The relative number of extreme vectors in each local set is reported in \autoref{fig:se}. First, observe that in general the number of extreme vectors is high for method \mL, lower for method \mM and lowest for method \mU. Second, method \mL is a rather special case; due to the generation method, all vectors are extreme except for $p=2$. Here rounding the vectors to integers results in non-extreme (unsupported) vectors. Next, considering methods \mM and \mU, an increasing number of objectives leads to more extreme vectors, which may be the result of an increasing volume of the hypercube. Finally, for a fixed number of objectives, the relative number of extreme vectors decreases for larger local set cardinalities (the lines for \mM and \mU move south). 

An instance of the \msp is generated given $p$ objectives and $S$ local sets of the same cardinality. The local sets are chosen given a \emph{method configuration} where the same method is used for all local sets or a combination of methods \mL and \mU (\mL{}\mU). In this case, half of the local sets use the same method if $S$ is even and if $S$ is odd, then the method that is applied to one extra local set is chosen randomly. If $k$ local sets with the same configuration are needed in the MSP instance, then they are selected randomly among the ten local set instances generated.  

Five instances for each objective $p=2,\ldots, 5$, number of local sets $S = 2, \ldots 5$, cardinality $|\Yn^s| = 50, 100, 200, 300$ and configuration \mU, \mL, \mM and \mL{}\mU are generated, resulting in 1280 MSP instances in total.  All instances are available at the Multi-Objective Optimization Repository \citep{MOrepo, MOrepo-Lyngesen24}. Note that since the vectors for each local set and objective are in the interval $[0, 10000]$, the hypercube containing $\Yn$ is $[0, 10000S]^p$. That is, the volume increases as $S$ and $p$ increases.

\subsection{\rqF}\label{sec:emp1}

\begin{figure}
    \centering
    \includegraphics[width=0.99\linewidth]{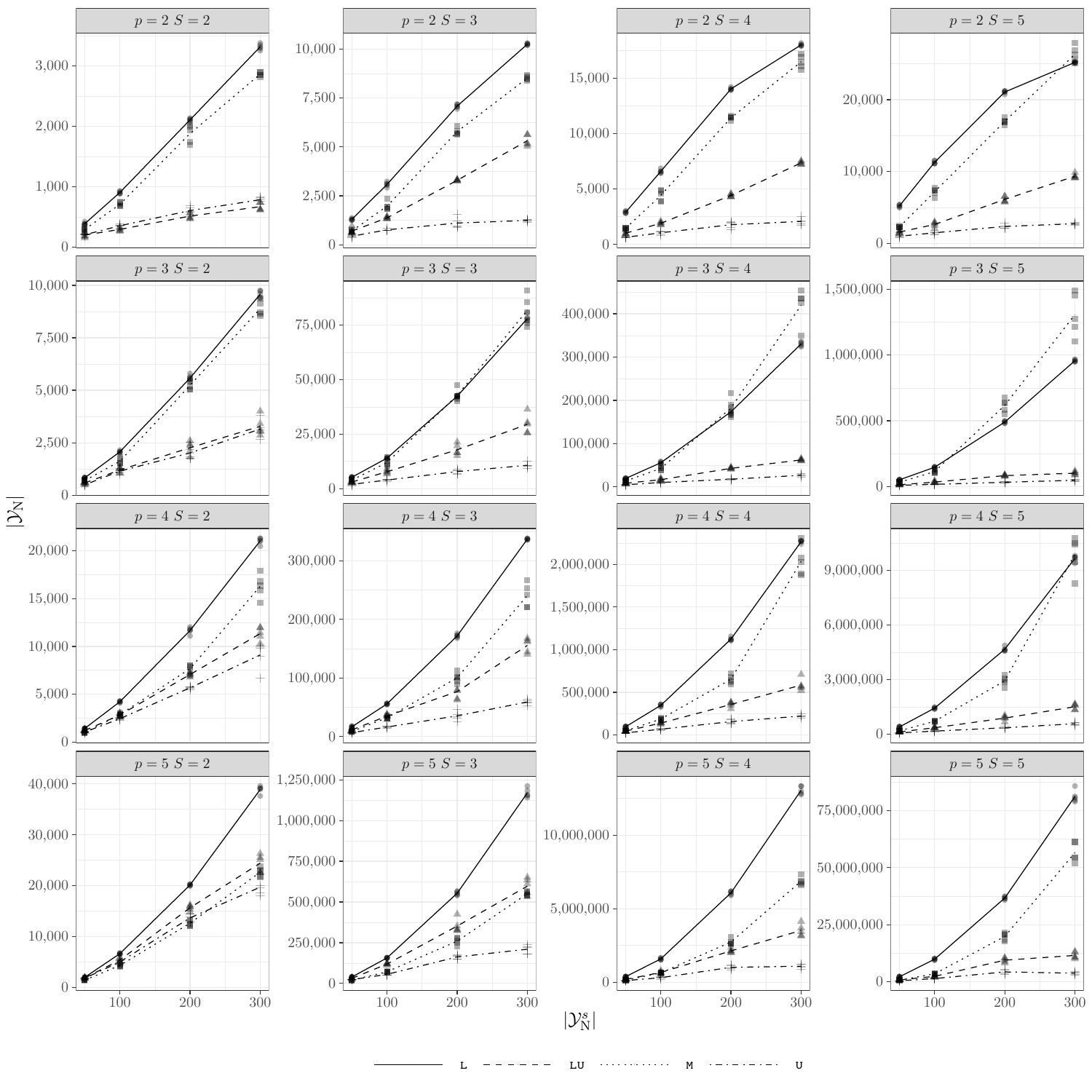}
    \caption{The cardinality of $\Yn$ given local set cardinality and fixed configuration. Averages over the 5 instances for each configuration are illustrated as lines. Axis ranges vary for each subplot.}
    \label{fig:Yn-size}
\end{figure}
In total all 1280 instances were solved. The cardinality of $\Yn$ is given in \autoref{fig:Yn-size}.

First, note that the cardinality of $\Yn$ grows as a function of local set cardinality, number of local sets, and number of objectives. Indeed, the \nd set becomes very large in some instances. Note that the volume of the hypercube containing $\Yn$ is $(10000|S|)^p$, due to the way the \nd set is generated for each local set, which may indicate an exponential relationship with $p$. To investigate, different functions were fitted to the results. A good fit was obtained using linear regression with function $|\Yn| = c_1 |\Yn^s|^{c_2p} S^{c_3p}$ (log transformed) given each configuration (R squared in range $[0.8, 0,95]$,  \cite{MOrepo-Lyngesen24}). That is, cardinality grows exponentially with $p$, and polynomial with $|\Yn^s|$ and $S$ with the highest effect for $S$ (the analysis resulted in $c_3 > c_2$).

Second, observe that the cardinality of $\Yn$ is affected by the shape of $\Yn^s$. The largest \nd sum is found when taking the MS of local sets with many extreme vectors (\mL) and lowest when the local sets have few extreme vectors (\mU). Furthermore in general, by comparing method \mM and \mLU it can be seen that finding the \nd sum using local sets with many extreme vectors in all local sets results in a higher cardinality of $\Yn$ (\mM).

The number of extreme vectors was found using an ``in hull'' algorithm \citep{gMOIP}. Due to memory issues, extreme vectors were found only for in total 1162/1280 instances. The relative number of extreme vectors in $\Yn$ was on average 4.5\% with a range of $[0.01\%, 33\%]$, decreasing with increasing $S$ and local set cardinality. This may be due to an increasing cardinality of the \nd sum. A decreasing relative number of extreme vectors for increasing problem sizes has also been observed in \cite{Sayin24}.

The highest average was obtained for method \mL (7.6\%) and lowest for method \mU (1.3\%). These rather low percentages imply that even though the relative number of extreme local vectors is close to 100\% (method \mL), this is not the case for the \nd sum, which contains many non-extreme vectors. That is, the VS of extreme local vectors obtained using different search directions may result in a non-extreme \nd VS, as pointed out in the comment following \autoref{prop:se}.

\subsection{\rqS}\label{subsection:empirical_MGS}

In total, $1270$ (out of 1280) instances were solved. Of these only $64$ instances did not satisfy the condition on \autoref{alg:MGS-if-fixed-is-reduced} in \autoref{alg:MGS} and of these only four did not satisfy the generating property on \autoref{l:end} in \autoref{alg:MGS}. For the remaining instances, the IP problem in \ref{app:IP-MGS} was solved (using the GLPK solver) and all optimal solutions were found. 
Of the 1270 instances solved, 1267 had a unique minimum generator set. This indicates that it is rare for an \msp{s} to have different \mgs{s}. For the remaining 3 instances, each had two \mgs{s}.

For each instance and local set, the relative \mgs cardinality was found, and the \emph{average relative \mgs cardinality} was calculated for methods $m\in\{\mU, \mL, \mM\}$:
\begin{equation}
  r^m = 100\frac{1}{|\S^m|}\sum_{s\in\S^m} \frac{|\G^s|}{|\Y^s|},
\end{equation}  
where $\S^m = \{ s\in\S \mid \text{method $m$ is used} \}$. That is, for an instance with configuration \mU, $r^\mU$ denotes average percentage reduction of all generators (since all local sets use method \mU). However, for the configuration \mLU, $r^\mL$ and $r^\mU$ denote the average percentage reduction of the generators for local sets using methods \mL and \mU, respectively.

Taking into account all instances, the ranges (overall averages) of the average relative \mgs cardinality were $r^\mU \in [0.7\%,100\%]$ (70\%), $r^\mM\in[15\%,100\%]$ (73\%) and $r^\mL\in[95\%,100\%]$ (100\%). For example, on average only $70$\% of the \nd vectors in local sets generated using method \mU were needed to generate the \nd sum. Since method \mU generates local sets with many unsupported vectors, one may conjecture that the relative \mgs cardinality is lower in local sets with many unsupported vectors. However, the range for \mU is $[0.7\%,100\%]$, which indicates that there are other factors affecting the numbers.

\begin{figure}[tb]
    \centering 
    \includegraphics[width=0.99\linewidth]{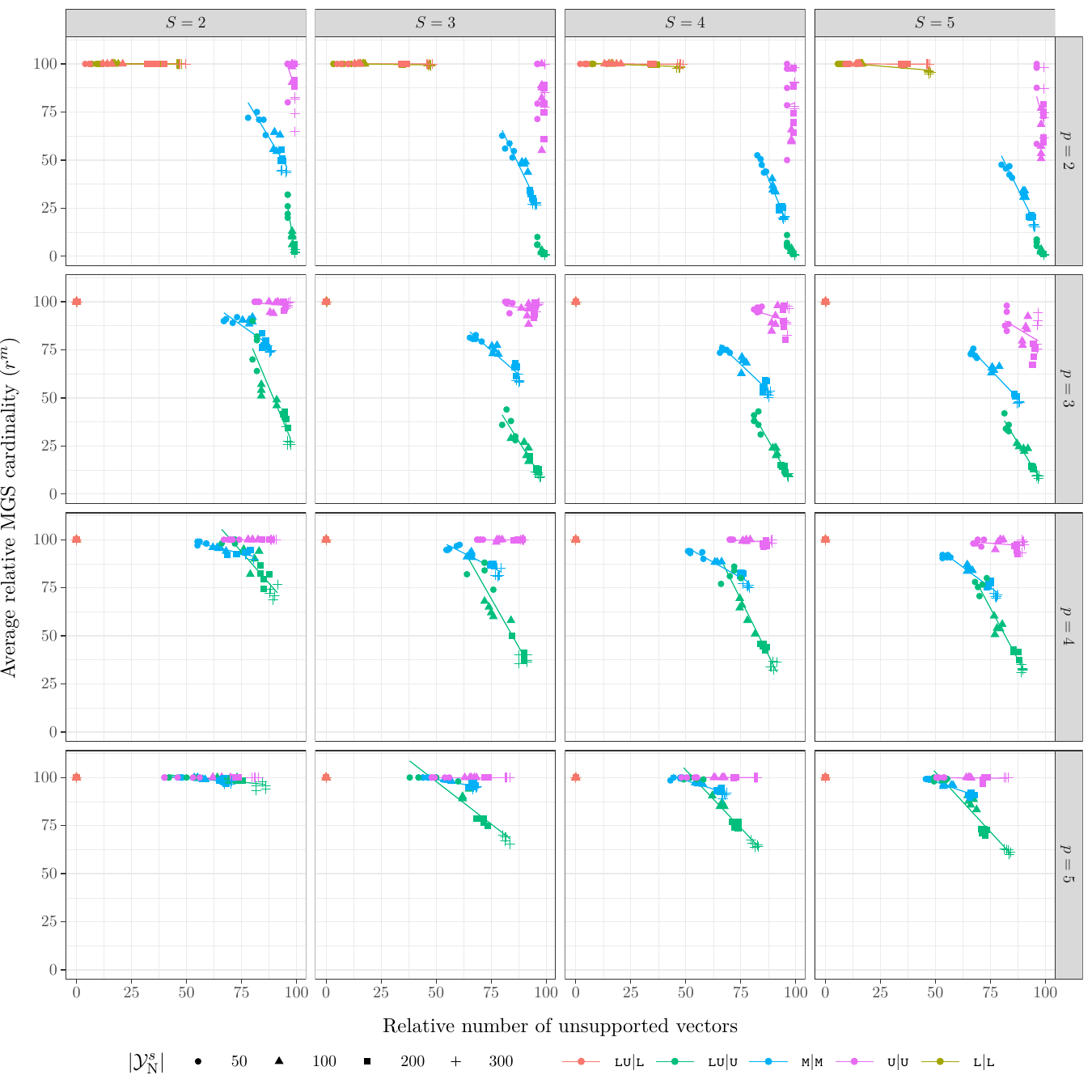}
    \caption{Relative \mgs cardinality given the relative number of unsupported \nd local vectors. Trend lines for each combination of configuration and local set method are given. Moreover, the cardinality of the local sets are visualized using different point shapes.}
    \label{fig:mgs-rel-values}
\end{figure}

To gain a deeper understanding of the results, consider \autoref{fig:mgs-rel-values}. Here the average relative \mgs cardinality $r^{m}$ for each instance is plotted against the relative number of unsupported \nd vectors in a local set. Colors are used to group results for instances generated using the same configuration and method. For example, \mLU{}$|$\mU corresponds to relative values for instances generated using configuration \mLU and local sets generated using method \mU. As a result, there are two vectors for each instance generated using configuration \mLU, since we have local sets generated using methods \mL and \mU. For all the other instances, we have only one point for each instance, since all local sets are generated using the same method.

First, observe that the \mLU{}$|$\mL and \mL{}$|$\mL cases are special, since local sets generated using method \mL, have almost only extreme vectors that are not redundant (see \autoref{fig:se}). Hence, due to \autoref{cor:extremeAreGenerators} the MGS is known in advance (generators equals $\Yn^s$ in most cases) and $r^m$ is approx. equal to 100. As a result, we do not include the \mLU{}$|$\mL and \mL{}$|$\mL cases in the analysis below.

Second, an increasing number of local sets results in a decrease in $r^m$ (a color group moves south across columns). That is, MSPs with more local sets give smaller MGSs relative to the cardinality of the local sets. The average given all instances of the relative \mgs cardinalities equals 81, 70, 67, 64 percent for $S$ equal 2, 3, 4 and 5, respectively. This is because more local sets increase the probability that a local vector is not part of a minimum generator set. For example, if $S=3$ then a \nd vector in local set 1 that is needed in the MS $(\Yn^1 \MinkSum \Yn^2)_\mcN = \Yn^{12}$, may be redundant in the MS $(\Yn^{12} \MinkSum \Yn^3)_\mcN = \Yn$. 

In contrast, there is a clear tendency that the number of local vectors required for a \mgs increases as a function of $p$ (a color group moves north between rows). The average given all instances of the relative \mgs cardinalities equals 65, 78, 90, 97 percent for $p$ equal 2, 3, 4 and 5, respectively. Recall that the cardinality of $\Yn$ grows exponentially with $p$ (see \autoref{sec:emp1}). Hence, more \nd local vectors may be needed for finding the \nd sum.

Also, observe that the relative \mgs cardinality is lower in local sets with many unsupported vectors if we consider each configuration separately. Indeed, the trend lines have a negative slope. However, the effect is very different. That is, the configuration has a large effect on $r^m$. For example, in the \mU{}$|$\mU case, much fewer \nd vectors are redundant compared to \mLU{}$|$\mU even though the number of unsupported vectors is the same. Furthermore, the \mM{}$|$\mM case has more redundant vectors compared to \mU{}$|$\mU even though this case has fewer unsupported vectors. This implies that the relationship between the local sets has a big impact on the number of redundant vectors. 

A deeper insight can be gained by considering \autoref{fig:ms-ex} which illustrates the MS of 4 instances ($p=2$ and $S=2)$. The plots also visualize the minimum generator set and the \nd sum. First, consider \autoref{fig:ms-uu}. Since both local sets are generated using the northeast part of a circle (\mU), in general, only a few vectors are redundant (on average 6\% in the instances). Next, compare this to \autoref{fig:ms-lu} where one local set is generated using the northeast part of a circle (\mU) and the second is generated using the southwest part of a circle (\mL). Here, some vectors are redundant for the \mU local set (on average 55\%). In between is the \mM instance given in \autoref{fig:ms-mm} where some vectors are redundant (on average 27\%). Finally, consider \autoref{fig:ms-lu-s}, similar to \autoref{fig:ms-lu}, except that one local set has been scaled so that the \nd set lies in a hypercube with different size. Notice that here only some of the vectors in the \mU local set are redundant compared to almost all in \autoref{fig:ms-lu}. That is, the combination of local set shapes has a big impact on the cardinality of the generator sets. %

\definecolor{Yone}{rgb}{0.8, 0.2, 0.2}
\definecolor{Ytwo}{rgb}{0.2, 0.6, 1}
\definecolor{Y}{rgb}{0.25098, 0.50196, 0}
\begin{figure}
\centering
\subfloat[MS (\mU and \mU).\label{fig:ms-uu}]{
	\includegraphics[width=0.43\linewidth]{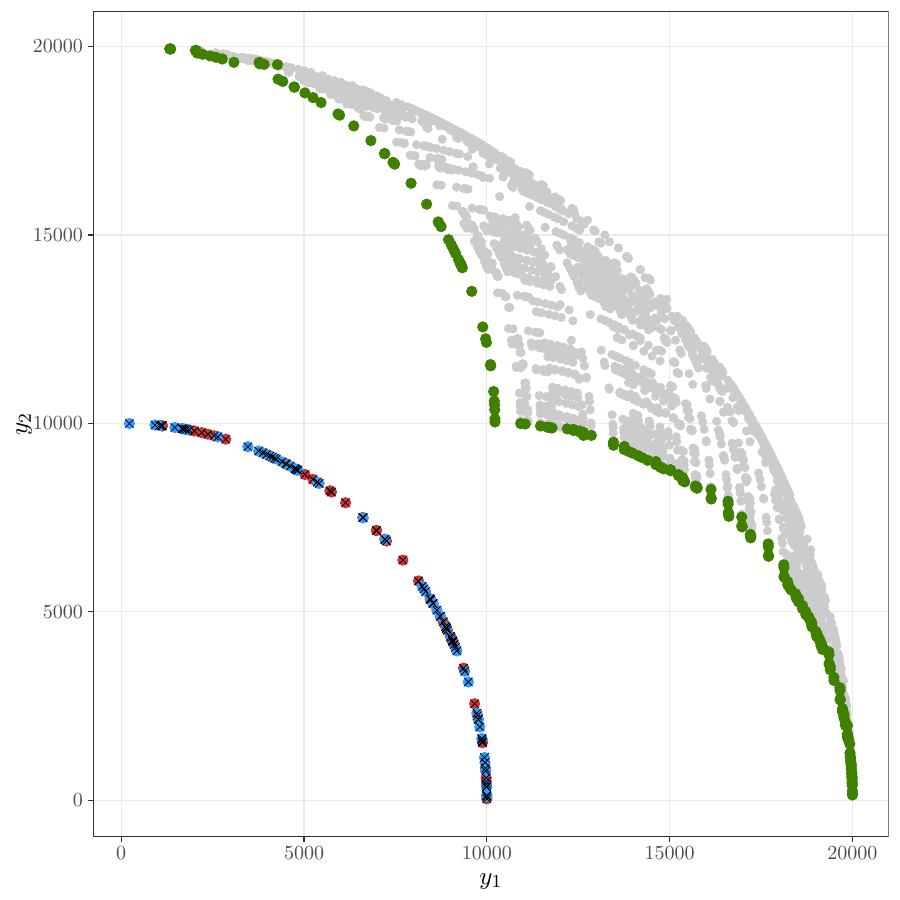}
} \hfill
\subfloat[MS (\mL and \mU).\label{fig:ms-lu}]{
	\includegraphics[width=0.43\linewidth]{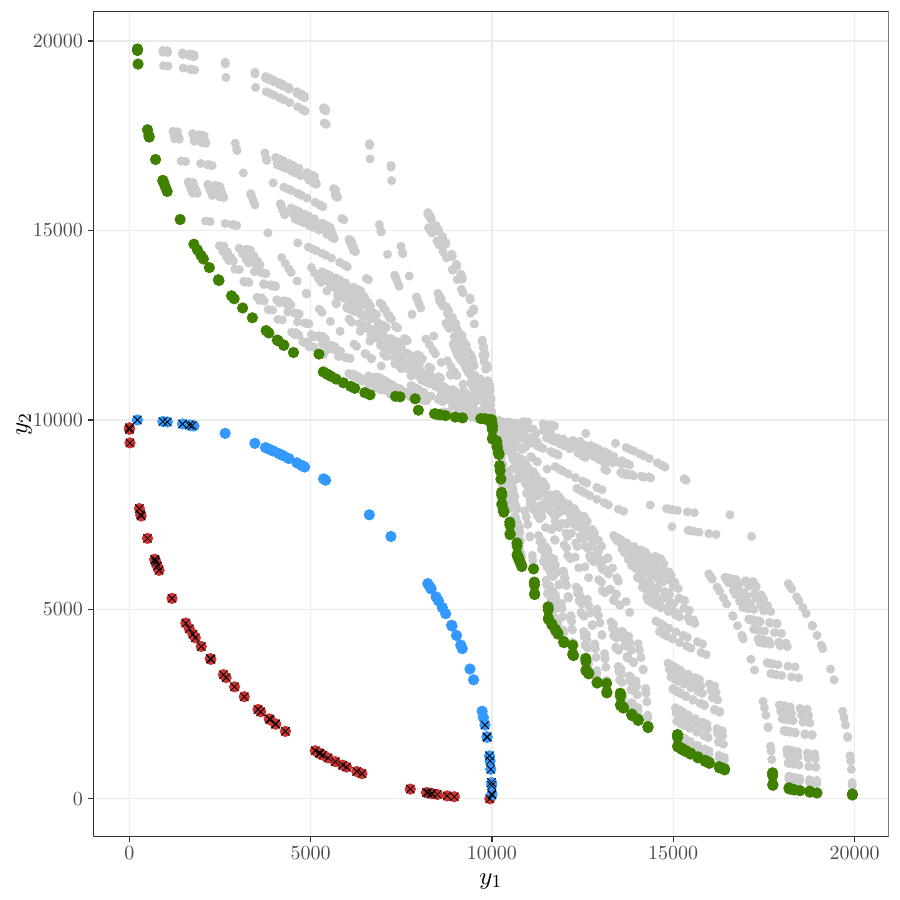}
} \\
\subfloat[MS (\mM and \mM).\label{fig:ms-mm}]{
	\includegraphics[width=0.43\linewidth]{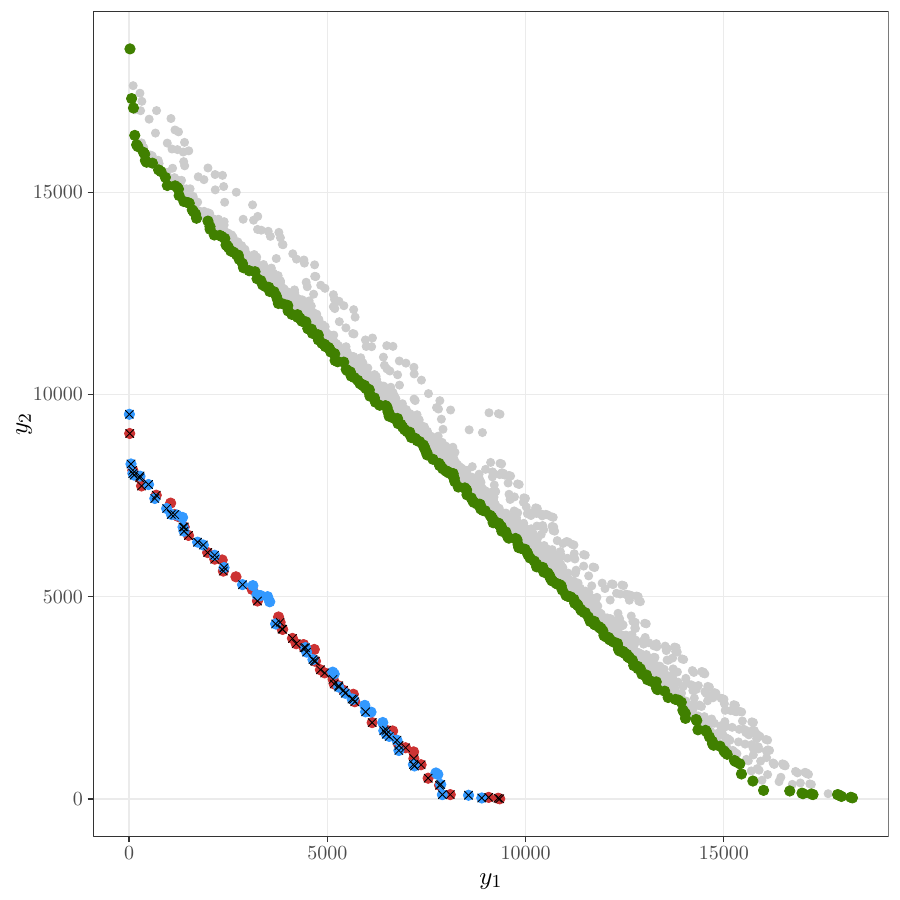}
} \hfill
\subfloat[MS (\mL and \mU, diffenent hypercubes).\label{fig:ms-lu-s}]{
	\includegraphics[width=0.43\linewidth]{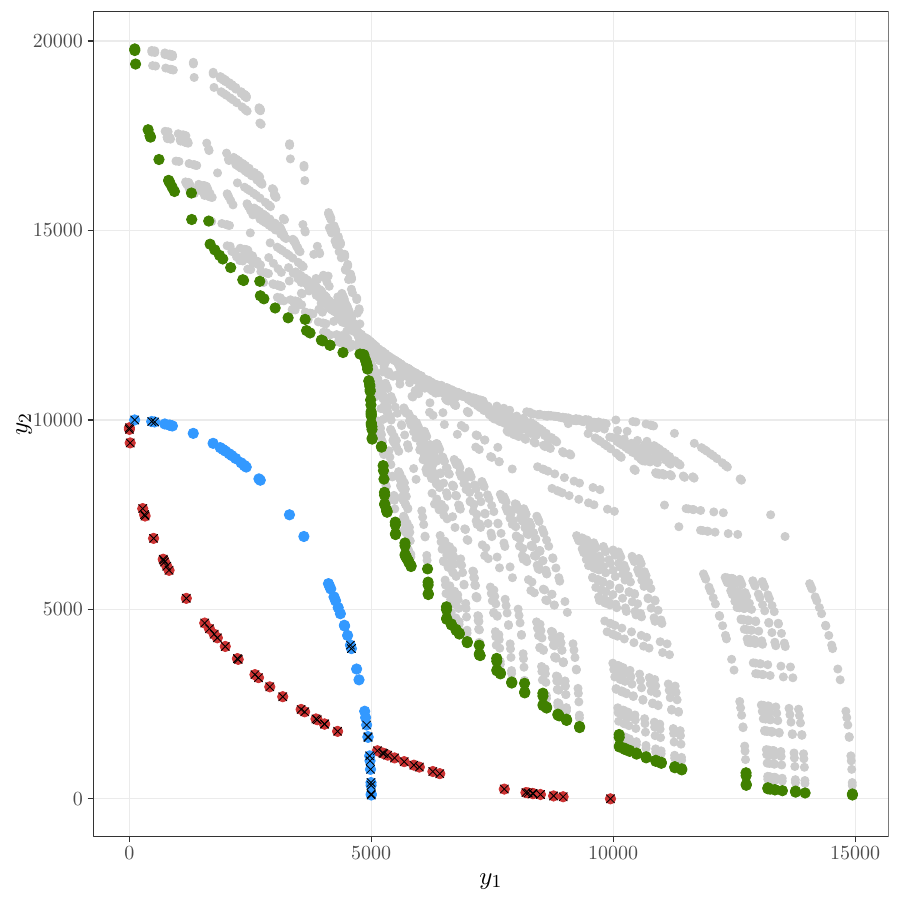}
}
\caption{Minkowski sum of two local sets \textcolor{Yone}{$\Yn^1$}, \textcolor{Ytwo}{$\Yn^2$} and the \nd sum \textcolor{Y}{$\Yn$}. The minimum generator set is visualized using a black color inside the local vectors.}\label{fig:ms-ex}
\end{figure}

\subsection{\rqT}

In this section we consider \autoref{alg:Alg3Pair} that apply pairwise comparisons between local sets to identify redundant vectors, and returns a generator set $\hat{\G}$. The algorithm assumes the \nd set for each local set $s$ is bounded by two stable sets $\L^s$ and $\U^s$ such that $\L^s \leqq \Yn^s \leqq \U^s$. That is, $\L^s$ and $\U^s$ are bounding the \nd set from \emph{below} and \emph{above}, respectively.

For measuring the quality of the set $\hat{\G}^s$ found using \autoref{alg:Alg3Pair} compared to a minimum generator $\G^s$ found in \autoref{subsection:empirical_MGS}, we consider the \emph{average relative number of redundant vectors} defined for each $m\in\{\mU, \mL, \mM\}$ as:
\begin{equation}\label{eq:rel-q}
  q^m = 100\frac{1}{|\S^m|}\sum_{\substack{s\in\S^m\\ |\Yn^s| \neq |\G^s|}} \frac{|\Yn^s|- |\hat \G^s|}{|\Yn^s|-|\G^s|},
\end{equation}  
where $\S^m = \{ s\in\S \mid \text{method $m$ is used} \}$. 
That is, if $q^\mU = 75$ for an instance then 75\% of all redundant vectors are found on average in the local sets generated using method \mU. Note that local sets with $\Yn^s = \G^s$ are excluded to get a fair comparison. One may argue that in this case the fraction should be one.

Since the largest MGS cardinality reductions are found for $p = 2$ objectives (see \autoref{subsection:empirical_MGS}), we limit the analysis to instances with $p = 2$ objectives and  $S = 2$ and $4$ local sets. That is, a total of 120 instances is considered. 

To test the hypothesis that the number of redundant vectors increases with the quality of the bounding sets, we define various bounding sets for each local set $s$. Since extreme vectors of a local set are never redundant, they are included in all bounding sets by assumption. In the following sections, we improve the bounding sets by considering $\L^s$ and $\U^s$, separately.

\subsubsection{Strengthening the bounding set from above}

\begin{figure}
    \centering
    \subfloat[$\L^s = \hull{\Ynse^s}_\textnormal{N}$ \label{fig:emp3-visual-a}]{
\includegraphics[width=0.48\linewidth]{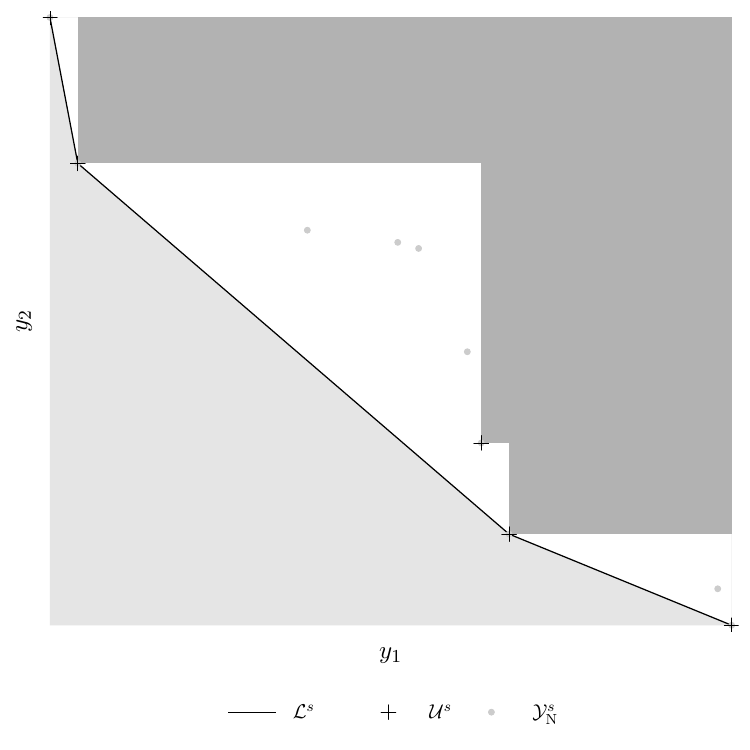}
    }
    \subfloat[$\L^s = \Yn^s$\label{fig:emp3-visual-b}]{
	\includegraphics[width=0.48\linewidth]{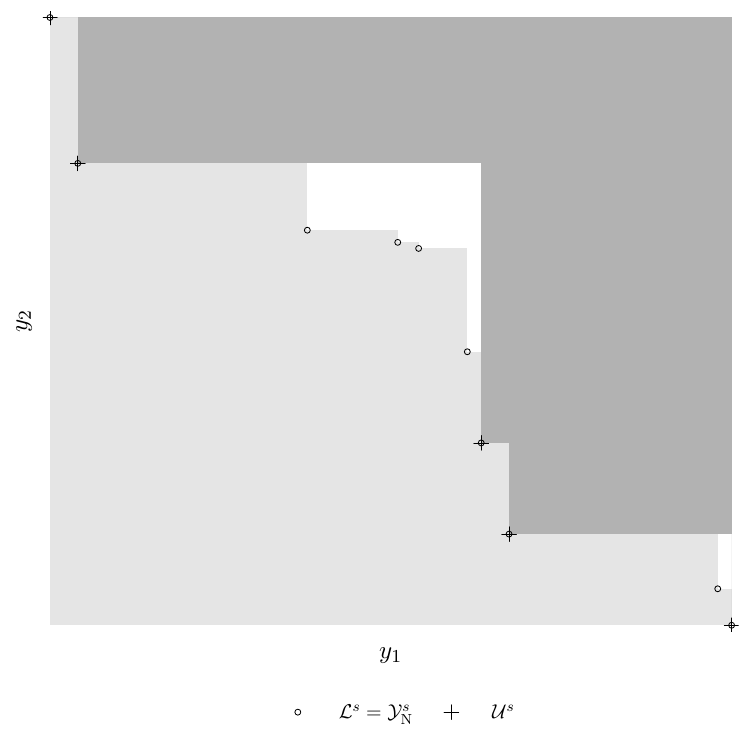}
    }
    \caption{Examples on different $\L^s$ bounding sets of a local set $s$. The light and dark gray area, excluding the boundary, illustrate where \nd vectors cannot be found, given $\L^s$ and $\U^s$, respectively.}
    \label{fig:emp3-visual}
\end{figure}

Assume that the bounding set $\L^s$ for each local set $s\in\S$ is fixed to $\L^s = \hull{\Ynse^s}_\textnormal{N}$ (see \autoref{fig:emp3-visual-a}). Moreover, let $\hat{\Y}^s(x)\subseteq \Yn^s \setminus \Ynse^s$ be a set consisting of $x$-percent of the non-extreme vectors in $\Yn^s$ (picked randomly).

To test \autoref{alg:Alg3Pair} on different $\U^s$, we introduce two parameters: given the local set $\bar s$ under consideration (line~\ref{alg:Alg3Pair-choose-s} in \autoref{alg:Alg3Pair}) let $\lambda \in [0,100]$ denote the percentage of non-extreme \nd vectors in $\U^{\bar s}$, \ie $\U^{\bar s} = \Ynse^{\bar s} \cup \hat\Y^{\bar s}(\lambda)$. Similarly, let $\gamma \in [0,100]$ denote the percentage of non-extreme \nd vectors in the other local sets, \ie $\U^{s} = \Ynse^{s} \cup \hat\Y^{s}(\gamma)$ for all $s\in\bar\S := \S \setminus \{s\}$ (line~\ref{alg:Alg3Pair-S-bar} in \autoref{alg:Alg3Pair}).
Hence, higher $\lambda$ values result in a stronger bound from above on the local set we seek to find redundant vectors in, and higher $\gamma$ values result in a stronger bound from above on the local sets we compare against (line~\ref{alg:Alg3Pair-if-dominated} in \autoref{alg:Alg3Pair}). %

\begin{figure}
    \centering
    \includegraphics[width=0.99\linewidth]{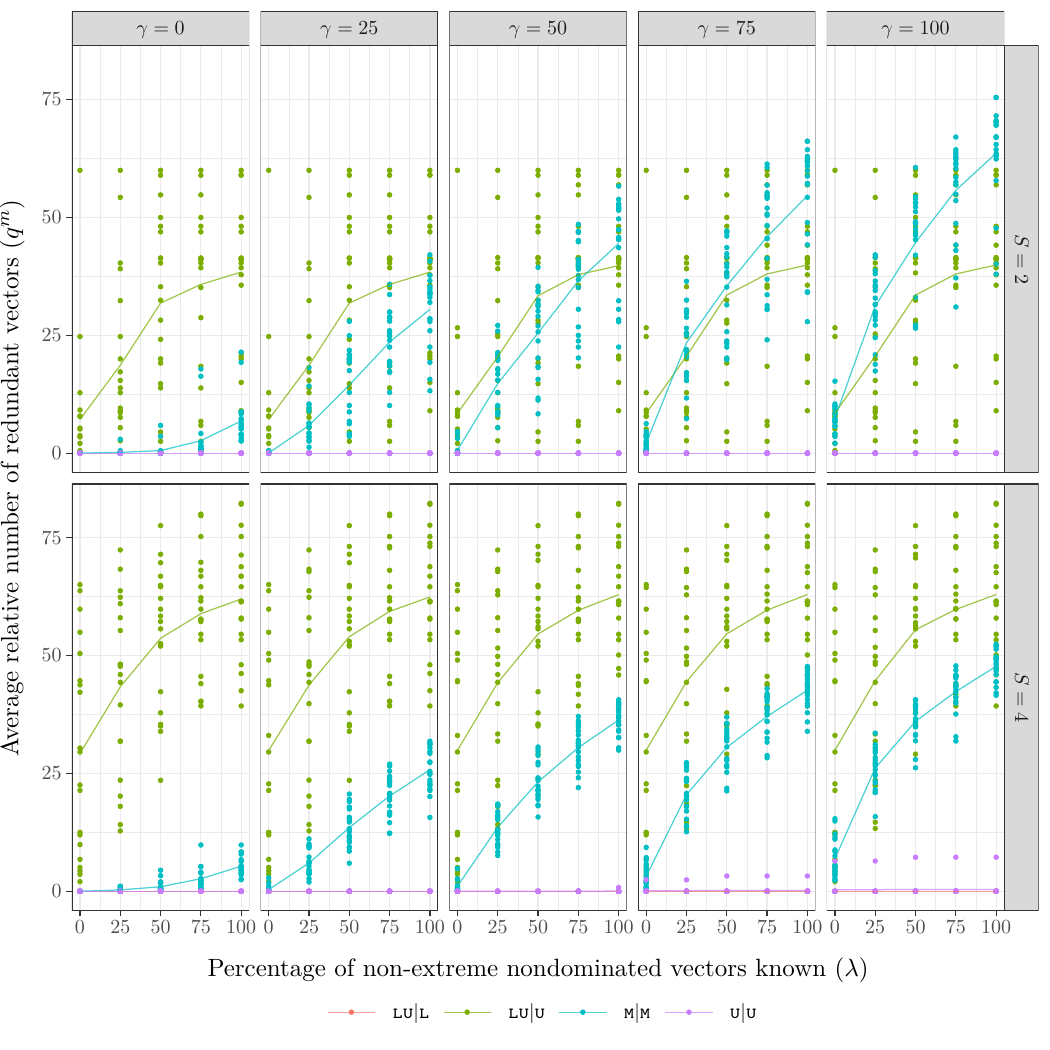}
    \caption{Comparing different bounding sets from above. Average relative number of redundant vectors ($q^m$) given percentage of non-extreme \nd vectors. Trend lines for each combination of instance configuration and local set method ($m$) are given.}
    \label{fig:emp3-plota}
\end{figure}

In our experiments, we calculate the generator set $\hat\G$ for different combinations of $\lambda$ and $\gamma$ in $\{0, 25, 50, 75, 100\}$.  %
The results are presented in \autoref{fig:emp3-plota}. 
Each subplot of \autoref{fig:emp3-plota} shows the average number of redundant vectors ($q^m$) for a fixed number of local sets ($S=2,4$) and a fixed value of $\gamma$ ($\gamma=0, 25, 50, 75, 100$).
As in \autoref{subsection:empirical_MGS}, colors are used to group results for instances generated using the same configuration and method. For example, \mLU{}$|$\mU corresponds to relative values for instances generated using configuration \mLU and local sets using method \mU.

First observe that for \mU{}$|$\mU and \mLU{}$|$\mL, the algorithm failed to identify any redundant vectors. Considering \mU{}$|$\mU local sets, $\Ynse^s$ often consists of only the lexicographic minima vectors. Hence, $\L^s$ is weak for all local sets $s\in\S$ and strengthening $\U^s$ has minimal effect. For \mLU{}$|$\mL local sets, \autoref{fig:mgs-rel-values} shows that only a few instances have $\Yn^s \neq \G^s$ (4 instances with 4 local sets). This indicates that redundant vectors are rare, and $q^\mL$  is undefined for most instances. Consequently, finding redundant vectors in this case is challenging.

Second, strengthening for the local set under consideration (increasing $\lambda$) leads to identifying more redundant vectors in the \mM{}$|$\mM and \mLU{}$|$\mU configurations, as indicated by the increasing lines in each subplot. Improved upper bounds therefore enhance the detection of redundant vectors. The effect varies; for example, for $S=4$
and $\lambda=\gamma=0.5$, the range (mean) of $q^\mM$ is $[8\%,36\%]$ ($23\%$).

Next, note that for \mM{}$|$\mM strengthening the other bounding sets (increasing $\gamma$) results in the identification of more redundant vectors. That is, having a better bound from above in the local set we compare against helps in identifying redundant vectors (the line moves north given the subplots in a row). 

This behavior does not apply to \mLU{}$|$\mU, as the lines are stable given the subplots in a row. For $S=2$, when comparing against a local set generated using method \mL, increasing $\gamma$ has little effect on strengthening $\U^{s}$ as $\vert \Ynse\vert \approx \vert \Yn\vert$. Thus, the behavior remain unchanged for varying $\gamma$. For $S = 4$ this also holds even though we also compare against a local set generated using method \mU (we compare against 2 generated using method \mL and one generated using method \mU). For this local set the bounding set from below is weak which makes it hard to find redundant vectors.

Finally, consider the effect of instances with more local sets (a fixed column in \autoref{fig:mgs-rel-values}). For the \mM{}$|$\mM case, we see a clear tendency that $q^\mM$ decreases as a function of $S$. The opposite is true for the \mLU{}$|$\mU case where $q^\mU$ increases as a function of $S$.
As such, the effect of $S$ is inconclusive.

\subsubsection{Strengthening the bounding set from below}

\begin{figure}
    \centering
    \includegraphics[width=0.99\linewidth]{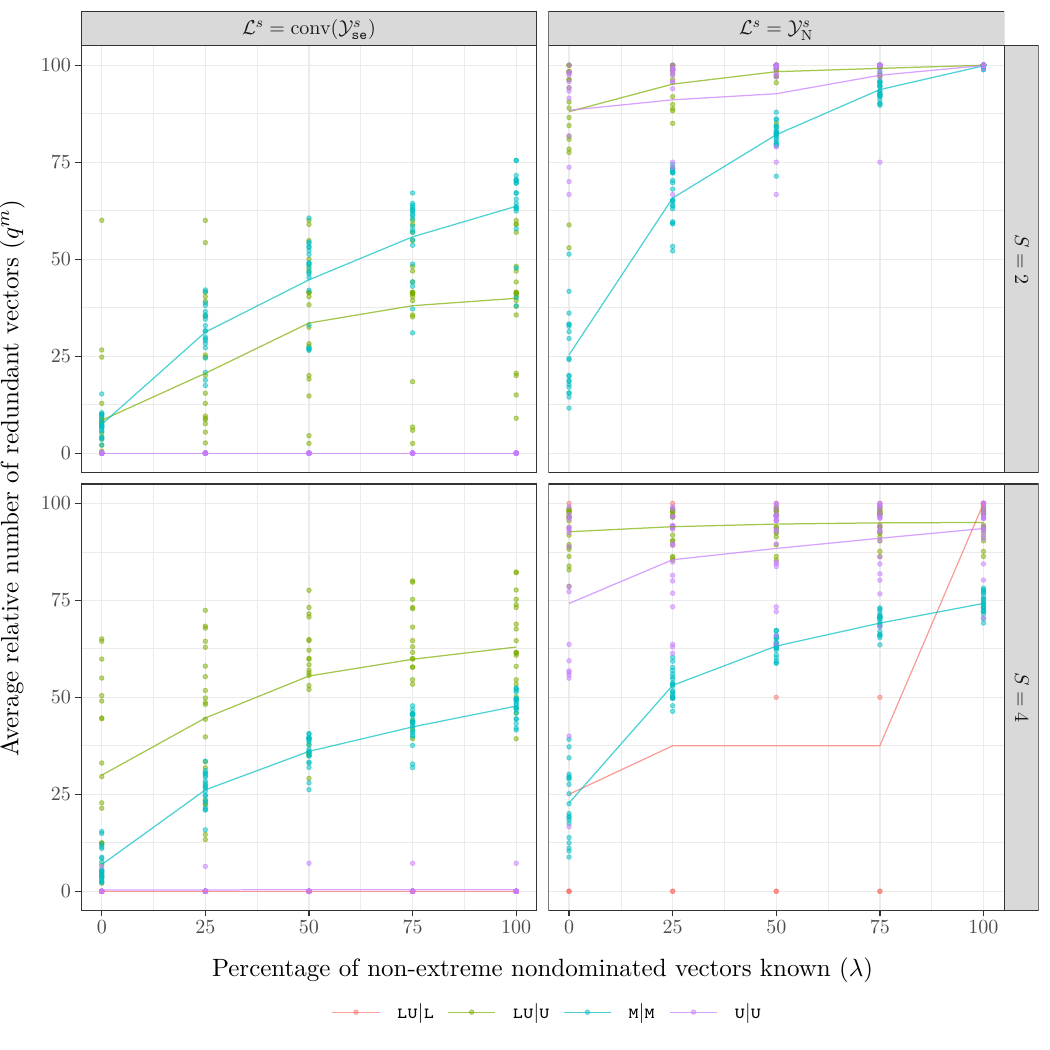}
    \caption{Comparing different bounding sets from below. Average relative number of redundant vectors ($q^m$) given fraction of non-extreme \nd vectors ($\lambda$, $\gamma = 1$). Trend lines for each combination of instance configuration and local set method ($m$) are given.}
    \label{fig:emp3-plotalb}
\end{figure}

The effect of strengthening the bounding sets from below is examined by considering the effect of changing the bounding sets $\L^{s}$ for each local set $s\in\bar\S$ compared against local set $\bar s$ from $\hull{\Ynse^{s}}_\textnormal{N}$ to $\Yn^{s}$ (line~\ref{alg:Alg3Pair-if-dominated} in \autoref{alg:Alg3Pair}). The two bounding sets are visualized in \autoref{fig:emp3-visual}. The results are given in \autoref{fig:emp3-plotalb} assuming the bounding set $\U^{s}$ for all $s\in\bar\S$ is fixed to $\U^{s} = \Yn^{s}$ ($\gamma = 1$).

Observe that strengthening $\L^{s}$ increases the $q^m$ values, as indicated by the northward shift of lines within each row. This shows that improving the bounding sets from below helps identify more redundant vectors. Moreover, there is a large improvement in $q^\mU$ when considering the \mU{}$|$\mU case. Recall that in this case $\Ynse^s$ often consists of only the lexicographic minima vectors. That is, $\L^{s} = \hull{\Ynse^{s}}_\textnormal{N}$ is weak. Hence, using strong bounding sets have a significant effect. The same is true for the few \mLU{}$|$\mL local sets where redundant vectors are identified due to the stronger bounding sets from below.

Finally, even using strong bounds may result in $q^m < 1$ ($\L^{s} = \Yn^{s}$, $\lambda^s = 1$, $S = 4$) . This may be due to the fact that \autoref{alg:Alg3Pair} only makes pairwise comparisons. Hence, redundant vectors, which can only be identified in the MS of more than two local sets, would not necessarily be identified by \autoref{alg:Alg3Pair} for $S=4$.

\section{Conclusion}\label{sec:conclu}

In this paper, we introduced the concept of generator sets for \MSP{s} (\msp{}s). While previous research focuses on efficiently filtering the Minkowski sum of several nondominated (\nd) local sets, we consider generator sets. That is, sufficient subsets of the local sets that generate the \nd sum.

First, we presented theoretical results showing that extreme supported local vectors are required for any generator set. In particular, any extreme supported vector is the sum of a unique set of extreme supported local vectors. Conversely, non-extreme local vectors can be redundant in generating the \nd sum. 

Second, based on these theoretical results we proposed an algorithm for efficiently finding minimum generator sets (\mgs{}s). The algorithm apply a simple pruning procedure fixing local vectors known to be in any generators.
The algorithm is based on an IP model with a preprocessing step fixing variables corresponding to local vectors known to be either redundant or necessary. Furthermore, we showed that bounding sets can be used to identify redundant local vectors and proposed an algorithm based on bounding sets, which identifies redundant local vectors using pairwise comparisons.

Through computational tests on diverse instances, we analyzed how the shape of local sets impacts the \nd sum. We found that the cardinality of the \nd sum grows exponentially with the number of objectives and polynomially with the cardinality and the number of local sets. Additionally, we demonstrated that, depending on the shape of the local sets, the number of objectives, and the number of subproblems, many vectors become redundant. The results showed that the combination of local set shapes significantly impacts the cardinality of the generator sets. Given a shape combination the relative \mgs cardinality was observed to decrease (linearly) with the number of unsupported local vectors. 

Finally, we used an algorithm based on bounding sets to identify redundant vectors. 
Having observed a large proportion of redundant vectors in the bi-objective case, we investigated how many of these the proposed algorithm could identify. The results showed that bounding sets can be used to find redundant vectors. Moreover, as expected, the relative number of redundant vectors identified increased as the bounding sets improved. In particular, both strong bounds from below and above are important for identifying redundant vectors.

A possible research direction is to incorporate this research in a bi-objective optimization algorithm for decomposable problems. Removing redundant vectors from the search space of subproblems could reduce the overall solution time.

\section*{Acknowledgements}

The numerical results presented in this work were obtained at the Centre for Scientific Computing, Aarhus University. Many thanks to Bruno Lang, University of Wuppertal, for providing the C implementation of the LimMem algorithm.

\appendix
\pagebreak
\section{IP formulation of minimum generator set problem}\label{app:IP-MGS}

Here an IP formulation of the \mgs problem, \eqref{prob:MGS}, is presented. 
We consider an MSP instance with sets $\Y^s$ for $s \in \S$, and \nd sum $\Yn := \nondom{\msum_{s\in \S} \Y^s}$.
Recall the definitions of $\C$ in \autoref{def:combinationsC} and the sets $\hat \Y^s$ and $\bar \Y^s$ for $s \in \S$ defined in \autoref{subsection:generator-sets}. 
Let $\bar{\mathcal{Y}}_{\textnormal{N}} = \Yn \cap \msum_{s \in \S} \bar{\mathcal{Y}}^s$ and define index set $\J$ for the non-generated vectors $y^j \in \Yn \setminus  \bar{\mathcal{Y}}_\textnormal{N} $ and index sets $\I^s$ for vectors $\tilde y^s_i \in \hat\Y^s \setminus \bar{\mathcal{Y}}^s$ for each $s \in \S$. With these index sets we consider the following decision variables:
\begin{align}
	x^s_i &= \begin{cases}
    	1, &\text{ if }\tilde{y}_i^s \in \G^s \text{ for } i \in \I^s\\ 
    	0, &\text{ otherwise }\\ 
    \end{cases}
    \\
		z^c_{j} &= \begin{cases}
    	1, &\text{ if $y^j \in \Yn$ is generated by combination $c \in \mathcal{C}(y^j)$} \\ 
    	0, &\text{ otherwise }\\ 
    \end{cases}
\end{align}

Letting $m_j^c\coloneqq\lvert \{s \in \S \mid c^s \in \hat{\Y^s} \setminus \bar{\Y^s}\}\rvert $ for each combination $c$ of \nd vector $y^j$, the \mgs problem, \eqref{prob:MGS}, can be formulated as the following integer linear programming problem:
\begin{align}
	\sum_{s \in \S}\lvert \bar{\Y^s}\rvert +\min\ &  \sum_{s \in \S}\sum_{i \in \I^s} x^s_i \\
	\text{s.t.:} \ & \sum_{\mathclap{c \in \mathcal{C}(y^j)}} z^c_j \ge 1 ,\quad &&\forall j \in \J \label{constr:generating}\tag{A.1}\\
	  \ & \sum_{s : c^s \in \hat{\Y^s} \setminus \bar{\Y^s} } \sum_{i \in \I^s}x^s_i \ge m_j^c z^c_j  && \forall c \in \mathcal{C}(y^j), \forall j \in \J\label{constr:generatingcomb}\tag{A.2}\\
		  \ & x^s_i \in \{0,1\} && \forall i \in \I^s, \forall s \in \S\\
		  \ & z^c_j \in \{0,1\} && \forall c \in \mathcal{C}(y^j), \forall j \in \J\\
\end{align}

The first set of constraints (\ref{constr:generating}) ensures that all vectors in the \nd sum are generated by at least one VS combination. The second set of constraints \eqref{constr:generatingcomb} ensures that a VS combination is valid only if each of the local vectors in the combination is chosen. A constant is added to the objective function (the first term) such that the objective value denotes the size of the \mgs.

Given an optimal solution $({\bar x},{\bar z})$ to the above IP a minimum generator set $\G = \{\G^1,\ldots,\G^S\}$ is given by $\G^s = \bar{\Y}^s \cup \{ y^s_i \in \Y^s\setminus \bar{\Y}^s : \bar{x}^s_i = 1\}$ for $s \in \S$.
Checking if the generator set is unique can be done by enumerating all optimal solutions to the above IP using, for example no-good inequalities.

\bibliographystyle{LaTeXConfig/elsarticle-harv} 
\bibliography{LaTeXConfig/literature}

\end{document}